# CARACTERIZACIÓN SEMÁNTICA DE LOS ÁRBOLES DE FORZAMIENTO SEMÁNTICO PARA LA LÓGICA DE PREDICADOS[a]

# SEMANTIC CHARACTERIZATION OF SEMANTIC FORCING TREES FOR PREDICATE LOGIC

MANUEL SIERRA-ARISTIZÁBAL[b*]

**RESUMEN:** La semántica de modelos para la lógica de predicados monádicos de primer orden, es caracterizada por una herramienta de inferencia visual llamada árboles de forzamiento semántico para la lógica de predicados monádicos. Las fórmulas que resultan válidas (o inválidas) mediante los árboles de forzamiento semántico, coinciden con las fórmulas válidas (o inválidas) mediante la semántica de valoraciones usual. En el caso que, una fórmula sea inválida mediante un árbol de forzamiento, un modelo que la refuta está determinado por las marcas de las hojas de este árbol. Este resultado es extendido a la semántica de modelos para la lógica de predicados diádicos de primer orden con dos variables.

**PALABRAS CLAVE:** Árbol de forzamiento semántico; predicados monádicos; semántica; valoración.

**ABSTRACT:** Model semantics for first-order monadic predicate logic is characterized by a visual inference tool called semantic forcing trees for monadic predicate logic. Formulas that are valid (or invalid) by semantic forcing trees match valid (or invalid) formulas by the usual valuation semantics. In the event that the formula is invalid by a forcing tree, a model that refutes it is determined by the marks of the leaves of this tree. This result is extended to model semantics for first-order dyadic predicate logic with two variables.

## 1. INTRODUCCIÓN

El método de las tablas semánticas es presentado por Beth (1962), y es sistematizado por Smullyan (1968), como árboles de opciones semánticas. El método explora sistemáticamente todas las posibilidades que podrían refutar una proposición dada, y determina cuál de ellas es lógicamente posible, en este caso se tiene un contraejemplo con el cual se refuta la validez de la proposición analizada. Si el contraejemplo no existe,

---

es decir ninguna posibilidad resulta viable, entonces la proposición analizada es válida. Este método tiene gran aceptación, como lo han hecho, entre otros, Carnielli (1987) para lógicas multivaluadas finitas, Barrero & Carnielli (2005) para la lógica clásica positiva, Areces *et al.* (2009) para memory logics, Cassano *et al.* (2015) para Calculus for Reasoning with Default Rules, Britz & VarzincZak (2019) para Contextual Defeasible ALC, Ferguson (2021) para weak Kleane logics, Bílková *et al.* (2021) para Two–Dimensional Fuzzy Logics, Indrzejczak & Zawidzki (2021) para Logics with Descriptions y Gr$a$̈tz (2021) para Non–deterministic Semantics.

Por otro lado, los árboles de forzamiento semántico presentados en Sierra (2001), no exploran todas las posibilidades en la búsqueda del contraejemplo, como se hace con los árboles de opciones, sino que, los árboles de forzamiento trabajan con las posibilidades que son deductivamente forzadas por un conjunto de reglas clara y completamente establecidas. Por esta razón, el análisis de validez con los árboles de forzamiento semántico es más simple y natural que el análisis con las tablas semánticas.

La prueba de la caracterización semántico–deductiva de los árboles de forzamiento, para el caso de la lógica proposicional, es presentada en Sierra (2006). Los árboles de forzamiento semántico para operaciones entre conjuntos son presentados en Sierra (2017) y los árboles de forzamiento semántico para la semántica de sociedades abiertas, correspondiente al sistema de lógica paraconsistente P1, son presentados en Sierra (2019).

En este trabajo, se demuestra detalladamente la caracterización semántica deductiva de los árboles de forzamiento, para el caso general de la lógica de predicados monádicos de primer orden (sin identidad ni funciones). En la sección 2 se presenta el lenguaje de la lógica de predicados monádicos de primer orden sin identidad ni funciones. En la sección 3 se presenta el algoritmo para construir el árbol inicial de una fórmula. En la sección 4 se presentan las reglas para marcar los nodos del árbol inicial, así como la generación de nuevas ramas; en esta sección también se define el concepto de validez desde el punto de vista de los árboles de forzamiento. En la sección 5 se presentan ejemplos representativos con el fin de ilustrar la aplicación de los árboles de forzamiento semántico, y la forma como se puede pasar de los árboles de un argumento válido, a la construcción de una prueba del argumento en el lenguaje natural; para el caso de árboles de argumentos inválidos, se muestra como a partir de las hojas del árbol se construye un modelo tradicional que refuta la validez del argumento. En la sección 6 se presentan los modelos para la lógica de predicados monádicos y diádicos junto con la definición de validez tradicional. En la proposición 5 de la sección 7, se prueba el objetivo central de este trabajo, es decir, la equivalencia de ambas nociones de validez.

## 2. LENGUAJE DE LA LÓGICA DE PREDICADOS

El lenguaje de la Lógica de Predicados, LP, consta de los conectivos binarios $\rightarrow$, $\wedge$, $\vee$ y $\leftrightarrow$, de los conectivos monádicos $\sim$, $\forall$ y $\exists$, además del paréntesis izquierdo y el paréntesis derecho. También se tiene una cantidad enumerable de variables, de constantes y de predicados $n$–ádicos (con $n = 1, 2, 3, \ldots$). El conjunto de fórmulas y de conjuntos de LP es generado por las siguientes reglas y sólo por ellas:

Si $z_1, \ldots, z_n$ son variables o constantes y $P$ es un predicado $n$–ádico entonces $P_{z_1, z_2, \ldots, z_n}$ es una fórmula atómica.

Toda fórmula atómica es una fórmula.

- Si $X$ es una fórmula entonces $\sim(X)$ es una fórmula.
- Si $X$ es una fórmula y x es una variable entonces $\forall(X)$ y $\exists(X)$ son fórmulas.
- Si $X$ y $Y$ son fórmulas entonces $(X) \wedge (Y)$, $(X) \vee (Y)$, $(X) \rightarrow (Y)$ y $(X) \leftrightarrow (Y)$ son fórmulas.

En las fórmulas $\forall x(X)$ y $\exists x(X)$, $\forall x$ y $\exists x$ son los cuantificadores universal y existencial para $x$, y $X$ es el alcance del cuantificador. Una ocurrencia de una variable $x$ se dice que es libre, si no se encuentra en el alcance de un cuantificador para $x$.

## 3. ÁRBOL DE UNA FÓRMULA

Sea $A$ una fórmula en la cual no figuran variables libres, el árbol inicial de A se representa por $Ar[A]$ y se construye utilizando las reglas que se presentan en la Figura 1 (donde $A$ y $B$ son fórmulas arbitrarias, $F(x)$ una fórmula con ocurrencias libres de la variable $x$, $F(\_)$ es el resultado de reemplazar las ocurrencias libres de $x$ en $F(x)$ por el espacio vacío "_"):

Por ejemplo, para el argumento: 'de $\forall(\sim P(x) \rightarrow Q(x))$ y $\exists x(P(x) \wedge Sa)$ se infiere $\exists x \, \forall y R(x, y)$ ', el condicional asociado es $[\forall x(\sim P(x) \rightarrow Q(x)) \wedge \exists x(P(x) \wedge Sa)] \rightarrow \exists x \forall y R(x, y)$ y su árbol inicial se presenta en la figura 2.

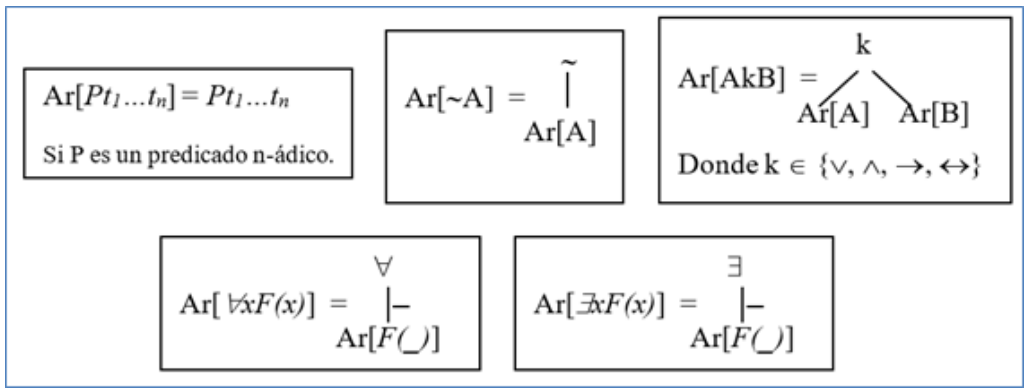

Figura 1: Reglas para la construcción del árbol inicial. Fuente: Elaboración Propia.

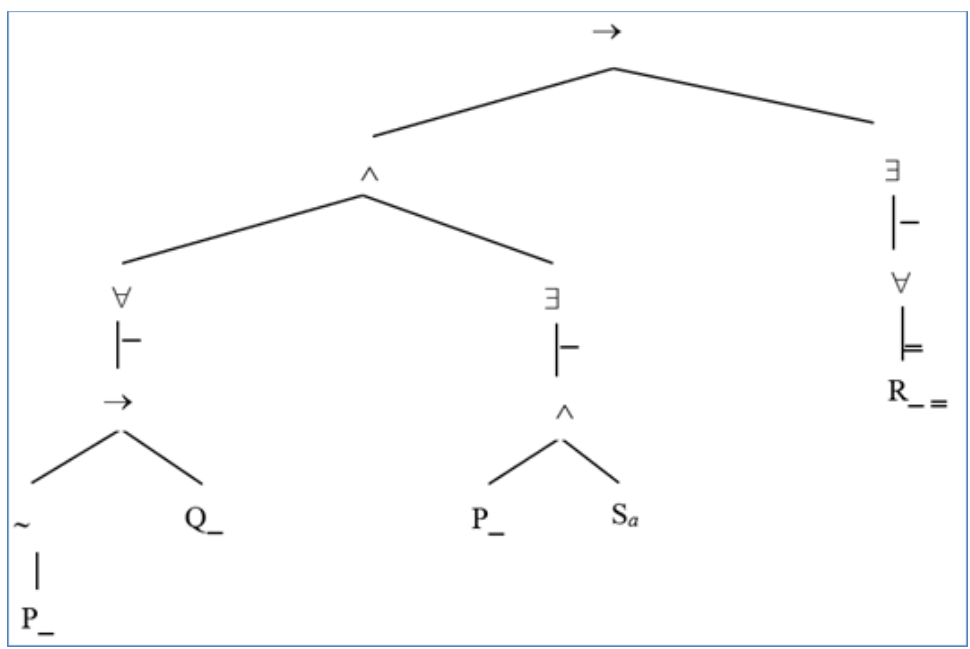

Figura 2: Árbol de la fórmula $(\forall x(\sim (Px \to Qx) \land \exists x(Px \land Sa)) \to \exists x \forall y Rxy$. Fuente: Elaboración Propia.

## 4. MARCANDO LOS NODOS DE UN ÁRBOL

Si un nodo $C$ es el conectivo monádico $\sim$, $\forall$ o $\exists$, entonces su único hijo se llama el alcance del operador y para hacer referencia a él se utiliza la notación $aC$.

Si un nodo $K$ es uno de los conectivos binarios $\land$, $\lor$, $\to$ o $\leftrightarrow$, entonces para sus hijos izquierdo y derecho se utiliza la notación $iK$ y $dK$, respectivamente.

Para toda fórmula $Y$, el nodo asociado a $Y$ es la raíz de $Y$, $R[Y]$, la cual a su vez es el operador principal de $Y$ en el caso que, $Y$ sea compuesta, o es la misma $Y$ en el caso que $Y$ sea atómica.

Para una fórmula $X$, $H(X)$ el conjunto de hojas del $Ar[X]$, y $N(X)$ el conjunto de nodos de $Ar[X]$.

Observación: En el caso de la fórmula $\exists x \forall y Rxy$, se deben asociar dos espacios vacíos, para diferenciar cual espacio vacío corresponde a la primera variable y cual espacio a la segunda, se diferencian los espacios mediante señales distintas, por ejemplo, para el primer espacio se utiliza la señal '_') y para el segundo la señal '$=$'). Por lo que, los espacios vacíos de $\exists x \forall y Rxy$, serían R $\_$ $=$ .

Para cada fórmula $X$, una función de marca de hojas, $m$ (o simplemente función de marca), es una función de $H(X)$ en $0, 1$.

- Si $m(p) = 1$ entonces se dice que la hoja $p$ está marcada con 1, o que es aceptada.

- Si $m(p) = 0$ entonces se dice que la hoja $p$ está marcada con 0, o que es rechazada.

Cada función de marca de hojas *m*, puede ser extendida de manera única (la prueba se presenta en la proposición 2), a una función de marca de nodos, *M*, de *N*(*X*) en 0, 1, haciendo *M*(*h*) = *m*(*h*) si *h* es una hoja, y aplicando las reglas de instanciación de espacios vacíos, junto con las reglas primitivas y derivadas para el forzamiento de marca, las cuales son presentadas a continuación en las secciones 4.1 a 4.5.

### 4.1. Reglas de instanciación o llenado de espacios vacíos

Los espacios vacíos del árbol inicial deben ser instanciados, es decir, deben ser llenados con constantes o con variables. Cuando un espacio vacío se instancia con una constante o con una variable, todos los espacios vacíos asociados a éste, se instancian con la misma constante o variable. Las instanciaciones de los espacios vacíos se realizan según los requerimientos dados por las siguientes reglas:

*IA*∀. Instanciación en la Aceptación del Universal: Si un universal es aceptado, entonces los espacios vacíos asociados a este universal (señales asociadas a la variable cuantificada), son instanciados para toda constante y/o para toda variable ya dadas, es decir, que ya figuran en el árbol que se está marcando. En el caso de no existir variables ni constantes debe ser introducida alguna de ellas, puesto que los modelos no pueden ser vacíos, es decir, se introduce una constante o se introduce una variable como representación de todas las constantes. En caso de realizarse más de una instanciación, entonces, del universal aceptado se deriva un hijo por cada instanciación.

*IR*∀. Instanciación en el Rechazo del Universal: Si un universal es rechazado, entonces los espacios vacíos asociados a este universal (señales asociadas a la variable cuantificada), son instanciados para alguna constante nueva (llamada testigo), es decir, una constante que no figura previamente en el árbol que se está marcando.

*I*∀. Instanciación en el Universal sin Marca: Si un universal no está marcado, entonces los espacios vacíos asociados a este universal (señales asociadas a la variable cuantificada), pueden ser instanciados para alguna constante o variable.

*IA*∃. Instanciación en la Aceptación del Existencial: Si un existencial es aceptado, entonces los espacios vacíos asociados a este existencial (señales asociadas a la variable cuantificada), son instanciados para alguna constante nueva (llamada testigo), es decir, una constante que no figura previamente en el árbol que se está marcando.

*IR*∃. Instanciación en el Rechazo del Existencial: Si un existencial es rechazado, entonces los espacios vacíos asociados a este existencial (señales asociadas a la variable cuantificada), son instanciados para toda constante y/o para toda variable ya dadas, es decir, que ya figuran en el árbol que se es-tá marcando. En caso de realizarse más de una instanciación, entonces, del existencial rechazado se deriva un hijo por cada instanciación.

*I∃.* Instanciación en el Existencial sin Marca: Si un existencial no está marcado, entonces los espacios vacíos asociados a este existencial (señales asociadas a la variable cuantificada), pueden ser instanciados para alguna constante o variable.

Observaciones:

- Las reglas *IR*∀ y *IA*∃ generan nuevas constantes. Esto se ilustra en los pasos 4 y 7 de la figura 6.

- Las reglas *IA*∀ y *IR*∃ aplican para las constantes nuevas y para las constantes del árbol inicial. Esto se muestra en los pasos 8 a 11 de la figura 6.

- Para el caso de las variables. La aplicación de las reglas *Aa*∀ o *Ra*∃ requiere variables en el espacio vacío (Ejemplo 1: en la figura 8 para llegar a la marca 1 en la raíz, se requiere la marca 1 en la rama derecha de la raíz, para lograrlo, la regla *Aa*∀ requiere que se instancie una variable, lo cual significa que en el paso 6 se debe instanciar la misma variable. Ejemplo 2: en la figura 10 para llegar a la marca 1 en la raíz, se requiere la marca 0 en el existencial de la rama derecha de la raíz, para lograrlo, la regla *Ra*∃ requiere que se instancie una variable, lo cual significa que en el paso 3 se debe instanciar la misma variable).

- Respecto a la aplicación de las reglas *I*∀ o *I*∃. Ejemplo 3: En la figura 10 para llegar al paso 6, la regla *Ra*∀ requiere que en el paso 4, a ser marcado con 1, se instancie una variable o constante que permita asegurar la marca 1, lo cual significa que en el paso 4 se debe instanciar la misma constante que fue instanciada en el paso 2. Ejemplo 4: En la figura 12 para llegar al paso 6, la regla *Aa*∃ requiere que en el paso 4, a ser marcado con 1, se instancie una variable o constante que permita asegurar la marca 1, lo cual significa que en el paso 4 se debe instanciar la misma constante que fue instanciada en el paso 3.

## 4.2. Reglas primitivas para el forzamiento de marcas de los nodos cuyos espacios vacíos ya han sido instanciados

*A∀.* Aceptación del Universal: Si un universal es aceptado, entonces el alcance es aceptado (para toda constante o variable previamente instanciada en el espacio vacío).

$$M(\forall) = 1 \Rightarrow M(a\forall) = 1.$$

Esta regla se ilustra en los pasos 10 y 11 de la figura 6. Observar que cuando se aplica esta regla a más de una variable o constante, el árbol inicial se modifica, puesto que de cada universal aceptado salen tantas ramas como variables y constantes haya, mientras que en el árbol inicial el universal sólo tiene un hijo.

*Aa*∃*.*     Aceptación del Alcance del Existencial: Si el alcance de un existencial es aceptado (para una constante o variable previamente instanciada en el espacio vacío), entonces el existencial es aceptado.

$$M(a\exists) = 1 \Rightarrow M(\exists) = 1.$$

Esta regla se ilustra en el paso 12 de la figura 7.

Al aplicar las reglas *OA* − *DM*, *OR* − *DM*, *OAi* − *Ad* →, *ORd* − *Ri* →, *ORi* − *Ad*∨, *ORd* − *Ai* ∨, presentadas más adelante en la sección 4.4, se tienen dos pasos de referencia, el inicial y el final. El paso inicial recibe el nombre de supuesto, y los pasos entre el inicial y el final, incluidos éstos, reciben el nombre de alcance del supuesto. Al aplicar las reglas mencionadas, se dice que el supuesto ha sido descargado. Cuando esto ocurre, las marcas de los pasos que están en el alcance del supuesto no pueden ser utilizadas en los pasos posteriores. Estas restricciones corresponden a las restricciones que se tienen en el sistema deductivo del cálculo de predicados, cuando se aplica el teorema de deducción o la prueba por contradicción, para los detalles ver Sierra (2006).

- Cuando $S(x)$ es un supuesto en el cual la variable $x$ ocurre libre, y $F(x)$ se encuentra en el alcance de $S(x)$, entonces se dice que $x$ en $F(x)$ no es independiente del supuesto $S(x)$, es decir, $x$ en $F(x)$ es dependiente del supuesto $S(x)$.

- Cuando $F(x)$ no se encuentra en el alcance de un supuesto, se dice que $x$ en $F(x)$ es independiente de tal supuesto.

- Cuando $S$ es un supuesto en el cual la variable $x$ no ocurre libre y $F(x)$ se encuentra en el alcance de $S(x)$, entonces $x$ en $F(x)$ es independiente del supuesto $S(x)$.

- Cuando la constante $a$ figura en $F(x)$, y la variable $x$ ocurre libre donde $a$ fue introducida (por aplicación de una de las reglas *IR*∀ y *R*∀ o *AA*∃ y *A*∃), entonces se dice que $x$ es dependiente del testigo $a$ en $F(x)$, en caso contrario, se dice que $x$ es independiente del testigo $a$ en $F(x)$.

- La variable $x$ es independiente en $F(x)$, si es independiente de todo testigo en $F(x)$, y si es independiente de todo supuesto en el cual $x$ ocurra libre.

*Aa*∀*.*     Aceptación del Alcance del Universal: Si el alcance de un universal es aceptado para una variable independiente en el espacio vacío, entonces el universal es aceptado.

$$M(a\forall) = 1 \Rightarrow M(\forall) = 1.$$

Esta regla se ilustra en el paso 13 de la figura 8, y también en el paso 5 de la figura 12 (observar la necesidad de que la variable a generalizar sea independiente).

$Ra\exists.$  Rechazo del Alcance del Existencial: Si el alcance de un existencial es rechazado para una variable independiente en el espacio vacío, entonces el existencial es rechazado.

$$M(a\exists) = 0 \Rightarrow M(\exists) = 0.$$

Esta regla se ilustra en el paso 7 de la figura 10.

### 4.3. Reglas derivadas para el forzamiento de marcas

Las reglas primitivas para el forzamiento de marcas son suficientes para estudiar las propiedades de los árboles de forzamiento, pero en la práctica, cuando se trata de marcar todos los nodos de un árbol, es importante tener reglas que cubran todas las posibilidades. A continuación, se presenta un juego completo de reglas derivadas (las reglas primitivas y derivadas para los conectivos $\wedge$, $\vee$, $\rightarrow$, $\leftrightarrow$ y $\sim$ se encuentran presentadas en Sierra (2006)).

Proposición 1. Reglas derivadas para los cuantificadores:

$R\forall.$  Rechazo del Universal: Si un universal es rechazado, entonces el alcance es rechazado (para una constante nueva previamente instanciada en el espacio vacío).

$$M(\forall) = 0 \Rightarrow M(a\forall) = 0.$$

Esta regla se ilustra en el paso 7 de la figura 6. Observar la necesidad de la constante nueva.

$Ra\forall.$  Rechazo del Alcance del Universal: Si el alcance de un universal es rechazado (para una constante o variable previamente instanciada en el espacio vacío), entonces el universal es rechazado.

$$M(a\forall) = 0 \Rightarrow M(\forall) = 0.$$

Esta regla se ilustra en el paso 6 de la figura 10.

$R\exists.$  Rechazo del Existencial: Si un existencial es rechazado, entonces el alcance es rechazado (para toda constante o variable previamente instanciada en el espacio vacío).

$$M(\exists) = 0 \Rightarrow M(a\exists) = 0.$$

Esta regla se ilustra en los pasos 8 y 9 de la figura 6. Observar que cuando se aplica esta regla a más de una variable o constante, el árbol inicial se modifica, puesto que de cada existencial rechazado salen tantas ramas como variables y constantes haya, mientras que en el árbol inicial el existencial sólo tiene un hijo.

$A\exists.$     Aceptación del Existencial: Si un existencial es aceptado, entonces el alcance es aceptado (para una constante nueva previa-mente instanciada en el espacio vacío).

$$M(\exists) = 1 \Rightarrow M(a\exists) = 1.$$

Esta regla se ilustra en el paso 10 de la figura 9. Observar la necesidad de la constante nueva.

Prueba de $R\forall$: Sea $M(\forall) = 0$. Supóngase que $M(a\forall) = 1$ para toda constante o variable (ya dadas) en el espacio vacío, entonces se tiene que $M(a\forall) = 1$ para una variable independiente en el espacio vacío, y por la regla $Aa\forall$ se infiere $M(\forall) = 1$, lo cual no es el caso.

Prueba de $Ra\forall$: Sea $M(a\forall) = 0$ para alguna variable o constante en el espacio vacío. Supóngase que $M(\forall) = 1$, entonces por la regla $A\forall$ se obtiene que $M(a\forall) = 1$ para toda variable o constante en el espacio vacío, lo cual contradice el supuesto inicial.

Prueba de $R\exists$: Sea $M(\exists) = 0$. Supóngase que $M(a\exists) = 1$ para alguna variable o constante en el espacio vacío, entonces por la regla $Aa\exists$ se infiere que $M(\exists) = 1$, lo cual no es el caso.

Prueba de $A\exists$: Sea $M(\exists) = 1$. Supóngase que $M(a\exists) = 0$ para toda constante o variable (ya dadas) en el espacio vacío, por lo que, $M(a\exists) = 0$ para una variable independiente en el espacio vacío, lo cual por la regla $Ra\exists$ significa que $M(\exists) = 0$, lo cual no es el caso.

### 4.4. Reglas para el forzamiento de marcas de los conectivos proposicionales

#### 4.4.1. Reglas primitivas

$A \wedge.$     Aceptación de la Conjunción: Si una conjunción es aceptada entonces tanto el hijo izquierdo como el derecho son aceptados.

$$M(\wedge) = 1 \Rightarrow [M(i\wedge) = 1 \text{ y } M(d\wedge) = 1].$$

*AiAd* $\wedge$ .      Aceptación a la Izquierda y Aceptación a la Derecha en la Conjunción: Si en una conjunción tanto el hijo izquierdo como el derecho son aceptados entonces la conjunción es aceptada.

$$[M(i\wedge) = 1 \text{ y } M(d\wedge) = 1] \Rightarrow M(\wedge) = 1.$$

*R* $\vee$ .      Rechazo de la Disyunción: Si en una disyunción es rechazada entonces tanto el hijo izquierdo como el derecho son rechazados.

$$M(\vee) = 0 \Rightarrow [M(i\vee) = 0 \text{ y } M(d\vee) = 0].$$

*RiRd* $\vee$ .      Rechazo a la Izquierda y Rechazo a la Derecha en la Disyunción: Si en una disyunción tanto el hijo izquierdo como el derecho son rechazados entonces la disyunción es rechazada.

$$[M(i\vee) = 0 \text{ y } M(d\vee) = 0] \Rightarrow M(\vee) = 0.$$

*R* $\rightarrow$ .      Rechazo del Condicional: Si un condicional es rechazado entonces el hijo izquierdo es aceptado y el hijo derecho es rechazado.

$$M(\rightarrow) = 0 \Rightarrow [M(i \rightarrow) = 1 \text{ y } M(d \rightarrow) = 0].$$

*AiRd* $\rightarrow$ .      Aceptación a la Izquierda y Rechazo a la Derecha en el Condicional: Si en un condicional el hijo izquierdo es aceptado y el hijo derecho es rechazado entonces el condicional es rechazado.

$$[M(i \rightarrow) = 1 \text{ y } M(d \rightarrow) = 0] \Rightarrow M(\rightarrow) = 0.$$

*A* $\leftrightarrow$ .      Aceptación del Bicondicional: Un bicondicional es aceptado si y solamente si ambos hijos tienen la misma marca, ambos son aceptados o ambos son rechazados.

$$M(\leftrightarrow) = 1 \Leftrightarrow M(i \leftrightarrow) = M(d \leftrightarrow).$$

*A* $\sim$ .      Aceptación de la Negación: Si una negación es aceptada entonces su alcance es rechazado.

$$M(\sim) = 1 \Rightarrow M(a \sim) = 0.$$

*Ra* $\sim$ .      Rechazo del Alcance de la Negación: Si el alcance una negación es rechazado entonces la negación es aceptada.

$$M(a \sim) = 0 \Rightarrow M(\sim) = 1.$$

### 4.4.2. Reglas derivadas

Las pruebas que garantizan la validez de estas reglas derivadas se encuentran en Sierra (2006).

*AiA* $\rightarrow$ .      Aceptación a la Izquierda y Aceptación del Condicional: Si son aceptados tanto el condicional como su hijo izquierdo entonces es aceptado el hijo derecho.

$$[M(i \rightarrow) = 1 \text{ y } M(\rightarrow) = 1] \Rightarrow M(d \rightarrow) = 1.$$

*RdA* $\rightarrow$ .      Rechazo a la Derecha y Aceptación del Condicional: Si un condicional es aceptado y su hijo derecho es rechazado entonces es rechazado el hijo izquierdo.

$$[M(d \rightarrow) = 0 \text{ y } M(\rightarrow) = 1] \Rightarrow M(i \rightarrow) = 0.$$

*Ri* $\rightarrow$ .      Rechazo a la Izquierda en el Condicional: Si en un condicional se rechaza el hijo izquierdo entonces se acepta el condicional.

$$M(i \rightarrow) = 0 \Rightarrow M(\rightarrow) = 1.$$

*Ad* $\rightarrow$ .      Aceptación a la Derecha en el Condicional: Si en un condicional se acepta el hijo derecho entonces se acepta el condicional.

$$M(d \rightarrow) = 1 \Rightarrow M(\rightarrow) = 1.$$

*AiR* $\wedge$ .      Aceptación a la Izquierda y Rechazo de la Conjunción: Si se rechaza una conjunción, pero se acepta su hijo izquierdo entonces se rechaza su hijo derecho.

$$[M(i\wedge) = 1 \text{ y } M(\wedge) = 0] \Rightarrow M(d\wedge) = 0.$$

$AdR \wedge$ .     Aceptación a la Derecha y Rechazo de la Conjunción: Si se rechaza una conjunción, pero se acepta su hijo derecho entonces se rechaza su hijo izquierdo.

$$[M(d\wedge) = 1 \text{ y } M(\wedge) = 0] \Rightarrow [M(i\wedge) = 0.$$

$Ri \wedge$ .     Rechazo a la Izquierda en la Conjunción: Si se rechaza el hijo izquierdo de una conjunción entonces se rechaza la conjunción.

$$M(i\wedge) = 0 \Rightarrow M(\wedge) = 0.$$

$Rd \wedge$ .     Rechazo a la Derecha en la Conjunción: Si se rechaza el hijo derecho de una conjunción entonces se rechaza la conjunción.

$$M(d\wedge) = 0 \Rightarrow M(\wedge) = 0.$$

$RiA \vee$ .     Rechazo a la Izquierda y Aceptación de la Disyunción: Si se acepta una disyunción, pero se rechaza su hijo izquierdo entonces se acepta su hijo derecho.

$$[M(i\vee) = 0 \, y \, M(\vee) = 1] \Rightarrow M(d\vee) = 1.$$

$Ai \vee$ .     Aceptación a la Izquierda en la Disyunción: Si se acepta el hijo izquierdo de una disyunción entonces se acepta la disyunción.

$$M(i\vee) = 1 \Rightarrow M(\vee) = 1.$$

$Ad \vee$ .     Aceptación a la Derecha en la Disyunción: Si se acepta el hijo derecho de una disyunción entonces se acepta la disyunción.

$$M(d\vee) = 1 \Rightarrow M(\vee) = 1.$$

*AiAd* $\leftrightarrow$ . Aceptación a la Izquierda y Aceptación a la Derecha en el Bicondicional: Si se aceptan ambos hijos de un bicondicional entonces se acepta el bicondicional.

$$[M(i \leftrightarrow) = 1 \, y \, M(d \leftrightarrow) = 1] \Rightarrow M(\leftrightarrow) = 1.$$

*RiRd* $\leftrightarrow$ . Rechazo a la Izquierda y Rechazo a la Derecha en el Bicondicional: Si se rechazan ambos hijos de un bicondicional entonces se acepta el bicondicional.

$$[M(i \leftrightarrow) = 0 \, y \, M(d \leftrightarrow) = 0] \Rightarrow M(\leftrightarrow) = 1.$$

*AiRd* $\leftrightarrow$ . Aceptación a la Izquierda y Rechazo a la Derecha en el Bicondicional: Si en un bicondicional se acepta el hijo izquierdo, pero se rechaza el hijo derecho entonces se rechaza el bicondicional.

$$[M(i \leftrightarrow) = 1 \, y \, M(d \leftrightarrow) = 0] \Rightarrow M(\leftrightarrow) = 0.$$

*AiA* $\leftrightarrow$ . Aceptación a la Izquierda y Aceptación del Bicondicional: Si se aceptan tanto el bicondicional como su hijo izquierdo entonces se acepta el hijo derecho.

$$[M(i \leftrightarrow) = 1 \, y \, M(\leftrightarrow) = 1] \Rightarrow M(d \leftrightarrow) = 1.$$

*RdA* $\leftrightarrow$ . Rechazo a la Derecha y Aceptación del Bicondicional: Si se acepta un bicondicional y se rechaza su hijo derecho entonces se rechaza el hijo izquierdo.

$$[M(d \leftrightarrow) = 0 \, y \, M(\leftrightarrow) = 1] \Rightarrow M(i \leftrightarrow) = 0.$$

*RiAd* $\leftrightarrow$ . Rechazo a la Izquierda y Aceptación a la Derecha en el Bicondicional: Si en un bicondicional se rechaza su hijo izquierdo y se acepta su hijo derecho entonces se rechaza el bicondicional.

$$[M(i \leftrightarrow) = 0 \, y \, M(d \leftrightarrow) = 1] \Rightarrow M(\leftrightarrow) = 0.$$

*RiA* $\leftrightarrow$ . Rechazo a la Izquierda y Aceptación del Bicondicional: Si se acepta un bicondicional y se rechaza su hijo izquierdo entonces se rechaza el hijo derecho.

$$[M(i \leftrightarrow) = 0 \, y \, M(\leftrightarrow) = 1] \Rightarrow M(d \leftrightarrow) = 0.$$

$AdA \leftrightarrow .$      Aceptación a la derecha y Aceptación del Bicondicional: Si se acepta un bicondicional y se acepta su hijo derecho entonces se acepta el hijo izquierdo.

$$[M(d \leftrightarrow) = 1 \, y \, M(\leftrightarrow) = 1] \Rightarrow M(i \leftrightarrow) = 1.$$

$RiR \leftrightarrow .$      Rechazo a la Izquierda y Rechazo del Bicondicional: Si se rechaza un bicondicional y se rechaza su hijo izquierdo entonces se acepta el hijo derecho.

$$[M(i \leftrightarrow) = 0 \, y \, M(\leftrightarrow) = 0] \Rightarrow M(d \leftrightarrow) = 1.$$

$AiR \leftrightarrow .$      Aceptación a la Izquierda y Rechazo del Bicondicional: Si se rechaza un bicondicional y se acepta su hijo izquierdo entonces se rechaza el hijo derecho.

$$[M(i \leftrightarrow) = 1 \, y \, M(\leftrightarrow) = 0] \Rightarrow M(d \leftrightarrow) = 0.$$

$RdR \leftrightarrow .$      Rechazo a la Derecha y Rechazo del Bicondicional: Si se rechaza un bicondicional y se rechaza su hijo derecho entonces se acepta el hijo izquierdo.

$$[M(d \leftrightarrow) = 0 \, y \, M(\leftrightarrow) = 0] \Rightarrow M(i \leftrightarrow) = 1.$$

$AdR \leftrightarrow .$      Aceptación a la Derecha y Rechazo del Bicondiciona: Si se rechaza un bicondicional y se acepta su hijo derecho entonces se rechaza el hijo izquierdo.

$$[M(d \leftrightarrow) = 1 \, y \, M(\leftrightarrow) = 0] \Rightarrow M(i \leftrightarrow) = 0.$$

$Aa \sim .$      Aceptación del Alcance de la Negación: Si el alcance de una negación es aceptado entonces la negación es rechazada.

$$M(a \sim) = 1 \Rightarrow M(\sim) = 0.$$

$R \sim .$      Rechazo de la Negación: Si la negación es rechazada entonces su alcance es aceptado.

$$M(\sim) = 0 \Rightarrow M(a \sim) = 1.$$

*IA.*      Iteración de la Aceptación: Sean *n* y *k* dos nodos asociados a una misma fórmula, si el nodo *n* es aceptado entonces el nodo *k* también es aceptado.

$$n \text{ asociado a } \beta, k \text{ asociado a } \beta, n \neq k \text{ y } M(n) = 1] \Rightarrow M(k) = 1.$$

*IR.*      Iteración del Rechazo: Sean *n* y *k* dos nodos asociados a una misma fórmula, si el nodo *n* es rechazado entonces el nodo *k* también es rechazado.

$$n \text{ asociado a } \beta, k \text{ asociado a } \beta, n \neq k \text{ y } M(n) = 1] \Rightarrow M(k) = 0.$$

*OA − DM.*      Opción de Aceptación que genera Doble Marca: Si al suponer que un nodo *N* está marcado con 1 (opción de aceptación *OA*) y al aplicar las reglas para marcar nodos, se tiene como consecuencia marcas diferentes en algún par de nodos asociados a una misma fórmula, entonces el nodo *N* realmente está marcado con 0.

Para cada nodo *n*,

$$[M(n) = 1 \Rightarrow \text{ para algún nodo } k, M(k) = 1 \text{ y } M(k) = 0] \Rightarrow M(n) = 0.$$

*OR − DM.*      Opción de Rechazo que genera Doble Marca: Si al suponer que un nodo *N* está marcado con 0 (opción de rechazo *OR*) y al aplicar las reglas para marcar nodos, se tiene como consecuencia marcas diferentes en algún par de nodos asociados a una misma fórmula, entonces el nodo *N* realmente está marcado con 1.

Para cada nodo *n*,

$$[M(n) = 0 \Rightarrow \text{ para algún nodo } k, M(k) = 1 \text{ y } M(k) = 0] \Rightarrow M(n) = 1.$$

Los pasos entre la opción de rechazo del nodo N y la contradicción, forman parte del llamado alcance del supuesto, y al aplicar la regla *OR − DM*, sólo se afirma aceptación del nodo *N*, y los otros pasos del alcance del supuesto son descargados, es decir, no pueden ser utilizados en pasos posteriores.

$RR - DM$. Rechazo de la Raíz que genera Doble Marca: Si en el árbol de la fórmula $\alpha$ se supone que la raíz está marcada con 0 (opción de rechazo $OR$) y al aplicar las reglas para marcar nodos, se tiene como consecuencia marcas diferentes en algún par de nodos asociados a una misma fórmula, entonces la fórmula $\alpha$ es $A$-válida. En este caso se dice que el árbol es un árbol mal marcado $AMM$.

Para $m$ una funcion de marca,

$$[M(R[\alpha]) = 0 \Rightarrow \text{ para algún nodo } k, M(k) = 1 \text{ y } M(k) = 0] \Rightarrow \alpha \text{ es } A - \text{Válida}.$$

Los pasos entre la opción de rechazo de la raíz y la contradicción forman parte del llamado alcance del supuesto, y al aplicar la regla $RR - DM$, sólo se afirma la aceptación de la raíz, lo cual significa que la fórmula $\alpha$ no puede ser refutada, es decir, la fórmula $\alpha$ es $A$-válida, y los otros pasos del alcance del supuesto son descargados, es decir, las marcas de estos nodos no son determinadas.

La regla $RR - DM$ es frecuentemente utilizada (es la versión, en los árboles de forzamiento semántico, de la prueba por contradicción), por lo que en vez de utilizar la opción de rechazo de la raíz $OR$ (raíz marcada con cuadro punteado), se utilizan las reglas Rechazo de la Raíz $RR$ y Doble Marca $DM$ como reglas primitivas.

La regla $DM$ se ilustra en el paso 17 de la figura 7.

OAi-Ad$\rightarrow$. Opción de Aceptación a la Izquierda que genera Aceptación a la Derecha en un Condicional: Si se supone que el antecedente es aceptado (opción de aceptación del antecedente $OA$) y al aplicar las reglas para marcar nodos, se tiene como consecuencia que el consecuente también es aceptado, entonces el condicional realmente es aceptado.

$$[M(i \rightarrow) = 1 \Rightarrow M(d \rightarrow) = 1] \Rightarrow M(\rightarrow) = 1.$$

Los pasos entre la opción de aceptación del antecedente y la aceptación del consecuente forman parte del llamado alcance del supuesto, y al aplicar la regla $OAi - Ad \rightarrow$, sólo se afirma la aceptación del condicional, y los pasos del alcance del supuesto son descargados, es decir, no pueden ser utilizados en pasos posteriores. Esta es la versión, en los árboles de forzamiento semántico, de la prueba condicional o teorema de deducción.

La regla $OAi - Ad \rightarrow$. se ilustra en el paso 14 de la figura 8.

ORd-Ri$\rightarrow$. Opción de Rechazo a la Derecha que genera Rechazo a la Izquierda en un Condicional: Si se supone que el consecuente es rechazado (opción de rechazo del consecuente $OR$) y al aplicar las reglas para marcar nodos, se tiene como consecuencia que el antecedente también es rechazado, entonces el condicional realmente es aceptado.

$$[M(d \to) = 0 \Rightarrow M(i \to) = 0] \Rightarrow M(\to) = 1.$$

Los pasos entre la opción de rechazo del consecuente y el rechazo del antecedente forman parte del llamado alcance del supuesto, y al aplicar la regla $ORd - Ri \to$, sólo se afirma la aceptación del condicional, y los pasos del alcance del supuesto son descargados, es decir, no pueden ser utilizados en pasos posteriores.

ORdi-Ad∨. Opción de Rechazo a la Izquierda que genera Aceptación a la Derecha en una Disyunción: Si se supone que el disyunto izquierdo es rechazado (opción de rechazo *OR*) y al aplicar las reglas para marcar nodos, se tiene como consecuencia que el disyunto derecho es aceptado, entonces la disyunción realmente es aceptada.

$$[M(i\lor) = 0 \Rightarrow M(d\lor) = 1] \Rightarrow M(\lor) = 1.$$

*ORd − Ai ∨ .* Opción de Rechazo a la Derecha que genera Aceptación a la Izquierda en una Disyunción: Si se supone que el disyunto derecho es rechazado y al aplicar las reglas para marcar nodos, se tiene como consecuencia que el disyunto izquierdo es aceptado, entonces la disyunción realmente es aceptada.

$$[M(d\lor) = 0 \Rightarrow M(i\lor) = 1] \Rightarrow M(\lor) = 1.$$

Los pasos entre la opción de rechazo de un disyunto y la aceptación del otro disyunto forman parte del llamado alcance del supuesto, y al aplicar la regla $ORi - Ad\lor$ o la regla $ORd - Ai$, sólo se afirma la aceptación de la disyunción, y los pasos del alcance del supuesto son descargados, es decir, no pueden ser utilizados en pasos posteriores.

### 4.5. *A*-Validez de una fórmula

Se dice que una fórmula *X* es *A*-válida (válida desde el punto de vista de los árboles) si y solamente sí para toda función de marca *m*, se tiene que $M(R[X]) = 1$.

Se dice que una fórmula *X* es *A*-inválida si no es *A*-válida, es decir si existe una función de marca *m*, tal que $M(R[X]) = 0$. En este caso, se dice que la función de marca refuta la fórmula *X*. También se dice que el árbol de *X* está bien marcado (*ABM*, todos sus nodos están marcados de acuerdo con las reglas sin generar contradicciones).

Observación. De lo anterior se puede afirmar que una fórmula $X$ es $A$-válida (válida desde el punto de vista de los árboles) si y solamente sí el árbol de $X$ no está bien marcado, es decir, $M(R[X]) = 0$, pero al marcar los nodos de acuerdo con las reglas, se genera una contradicción (es decir, dos nodos asociados a una misma fórmula figuran con marcas diferentes).

### 4.6. Presentación gráfica de las reglas para el forzamiento de marcas

Un nodo encerrado en un círculo indica que el nodo está marcado con 1 (es aceptado), un nodo encerrado en un cuadro indica que el nodo está marcado con 0 (es rechazado). La configuración inicial se muestra sobre la línea punteada y las marcas iniciales se presentan en color azul, la configuración que resulta cuando se aplica la regla, se presenta debajo de la línea punteada y la nueva marca generada por la regla se presenta en color rojo.

Los gráficos para el caso de los conectivos proposicionales, presentados más arriba, se encuentran en Sierra (2010). Para el caso de los cuantificadores, los gráficos correspondientes a las reglas se presentan en la figura 5, ilustraciones 17 y 18.

Observación. Las reglas para las opciones y reglas de opciones para el condicional y la disyunción, los círculos y cuadros punteados indican que se toma una opción (Opción Afirmativa, $OA$, para el caso del círculo punteado. Opción de Rechazo, $OR$, para el caso del cuadro punteado). Los arcos en las reglas para las opciones indican que el valor supuesto ($OA$ u $OR$) debe ser cambiado. Los arcos son opcionales y pueden ser omitidos.

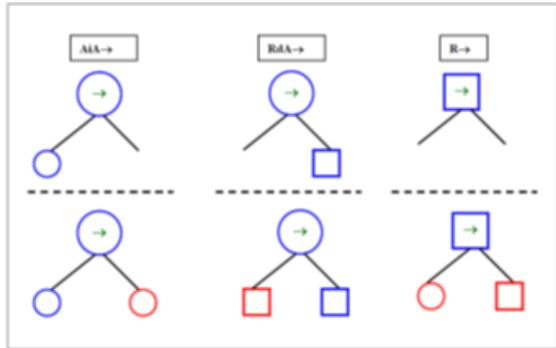

1. Reglas para el condicional. Fuente: Elaboración Propia.

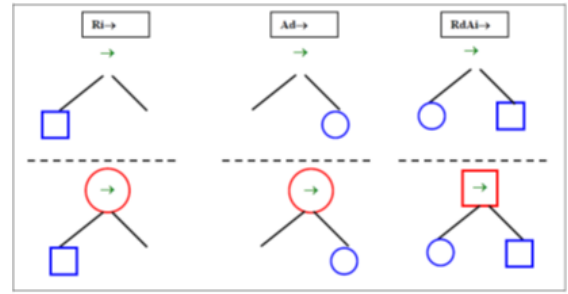

2. Reglas para el condicional. Fuente: Elaboración Propia.

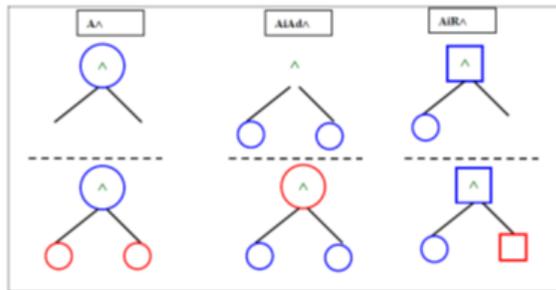

3. Reglas para la conjunción. Fuente: Elaboración Propia.

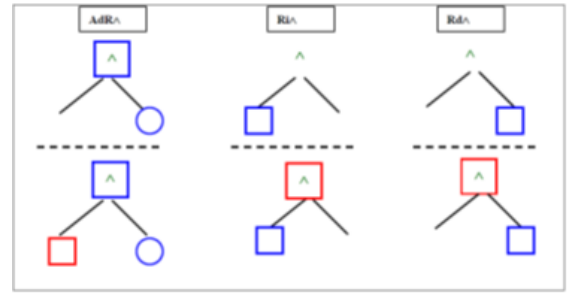

4. Reglas para la conjunción. Fuente: Elaboración Propia.

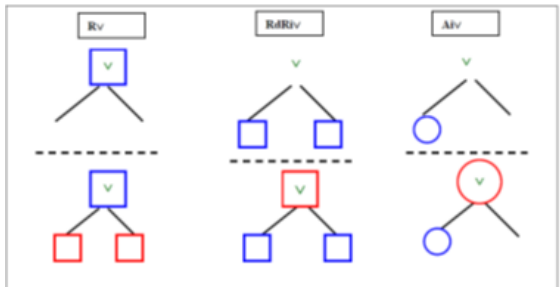

5. Reglas para la disyunción. Fuente: Elaboración Propia.

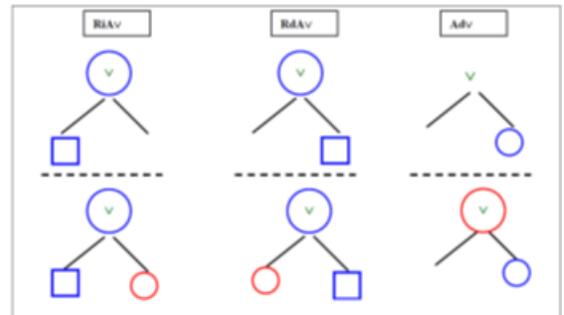

6. Reglas para la disyunción. Fuente: Elaboración Propia.

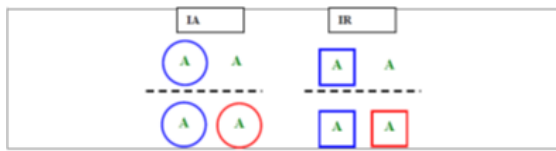

7. Reglas de iteración. Fuente: Elaboración Propia.

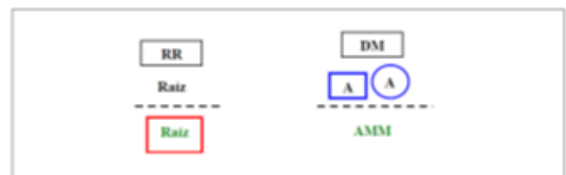

8. Reglas de iteración. Fuente: Elaboración Propia.

Figura 3: Reglas para los conectivos binarios, la interacción y la doble marca.

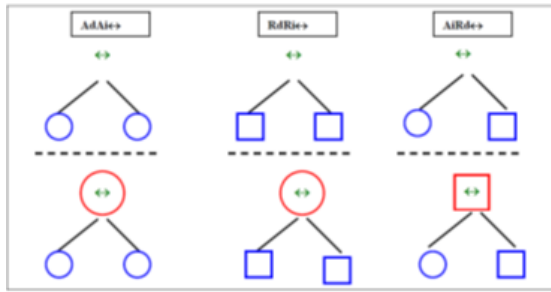

9. Reglas para el bicondicional. Fuente: Elaboración Propia

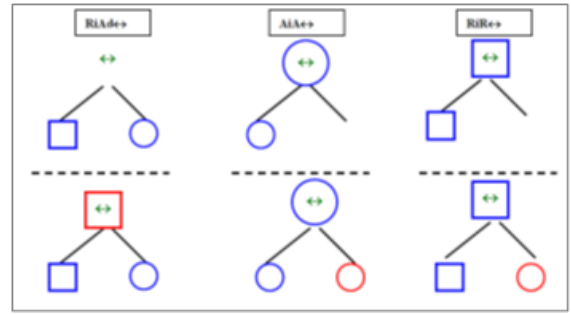

10. Reglas para el bicondicional. Fuente: Elaboración Propia

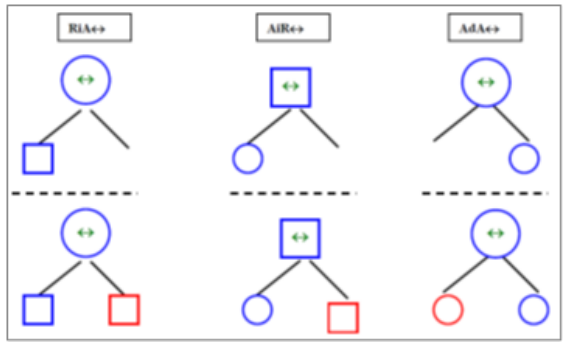

11. Reglas para el bicondicional. Fuente: Elaboración Propia

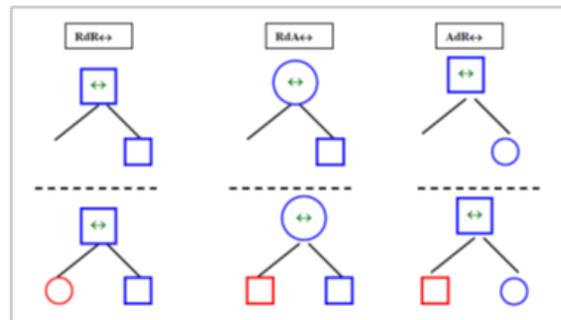

12. Reglas para el bicondicional. Fuente: Elaboración Propia

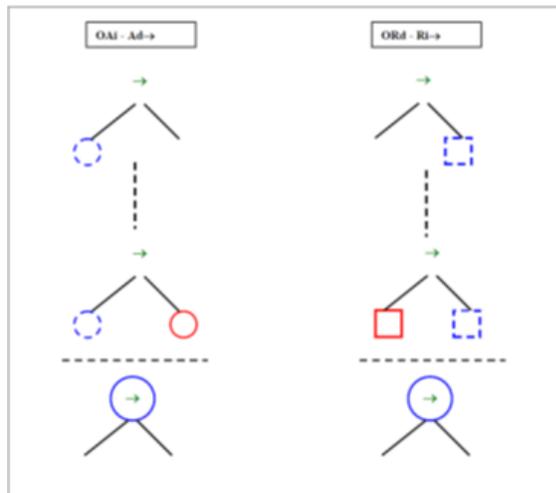

13. Reglas de opciones para el condicional. Fuente: Elaboración Propia

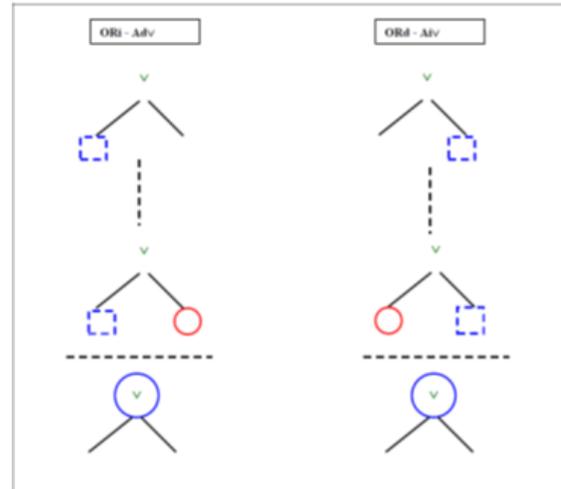

14. Reglas de opciones para la disyunción. Fuente: Elaboración Propia

Figura 4: Reglas para el bicondicional y reglas de opciones para el condicional y la disyunción.

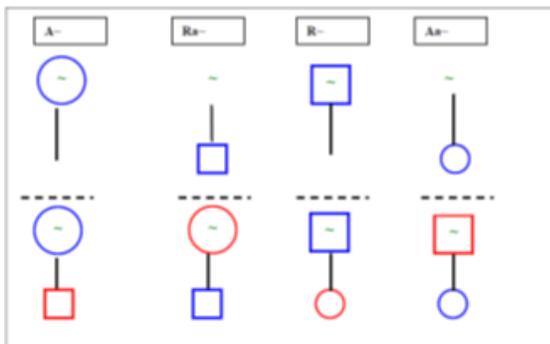

15. Reglas para la negación. Fuente: Elaboración Propia.

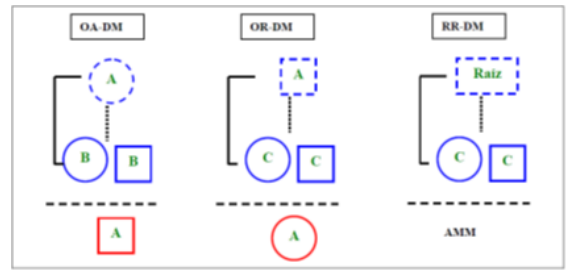

16. Reglas para cambio de opciones. Fuente: Elaboración Propia.

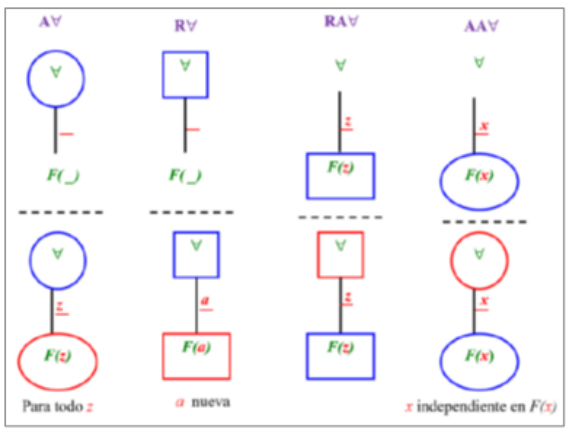

17. Reglas para el cuantificador universal. Fuente: Elaboración Propia.

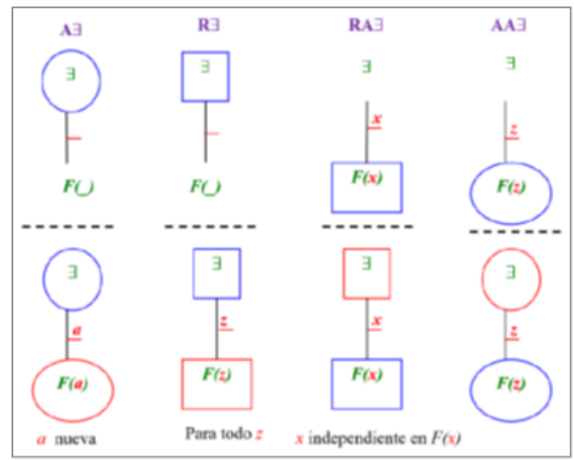

18. Reglas para el cuantificador existencial. Fuente: Elaboración Propia

Figura 5: Reglas para las opciones, para la negación y para los cuantificadores.

# 5. ILUSTRACIONES

Siguiendo a Sierra (2001), Sierra (2006), un árbol de forzamiento está mal marcado, cuando su raíz está marcada con 0, y además existen dos nodos asociados a una misma fórmula, los cuales tienen marcas contrarias. Si la raíz está marcada con 0, y todos los nodos están marcados sin generar contradicciones, se dice que el árbol está bien marcado, *ABM*. Lo anterior significa que cuando el árbol de una fórmula está mal marcado, entonces la fórmula es *A*-válida, y cuando el árbol de una fórmula está bien marcado, entonces la fórmula es *A*-inválida.

Observación: Cuando no se sabe si el argumento asociado al árbol es válido o inválido, se inicia el análisis utilizando la regla *RR* (rechazo de la raíz), y se determina si el árbol está bien marcado o no, lo cual permite concluir respectivamente si el argumento es inválido o no. Esta estrategia es conocida como forzamiento indirecto, y es utilizada en las ilustraciones 1, 2, 4 y 6.

Para efectos meramente argumentativos, cuando se sabe que el argumento es válido y la raíz del árbol es un condicional o una disyunción, además del forzamiento indirecto (el cual proporciona una prueba indirecta de la validez de la fórmula asociada al árbol), puede ser utilizado el llamado forzamiento directo (el cual proporciona una prueba directa de la validez de la fórmula asociada al árbol), el cual consiste en la aplicación de las reglas de toma de opciones OAi-Ad→, ORd-Ri→, ORi-Ad ∨ u ORd-Ai∨, con el objetivo de obtener la marca 1 en la raíz. Esta estrategia es utilizada efectivamente en las ilustraciones 3 y 5. En la ilustración 7 se intenta aplicar esta misma estrategia, pero según la ilustración 6, el argumento asociado al árbol es inválido, lo cual garantiza el fracaso de la estrategia de forzamiento directo en la ilustración 7, puesto que el forzamiento directo no constituye un método para determinar la validez de fórmulas.

## 5.1. Ilustración 1

En la figura 6 se muestra un árbol de forzamiento bien marcado para la fórmula *A*-inválida $\exists x(Px \land \forall yRxy) \to \forall x \exists yRxy$. Un nodo encerrado en un círculo indica que el nodo está marcado con 1 (es aceptado), un nodo encerrado en un cuadro indica que el nodo está marcado con 0 (es rechazado).

La fórmula *A*-inválida $\exists x(Px \land \forall yRxy) \to \forall x \exists yRxy$, corresponde a la formalización del siguiente argumento: Si el último trabajo de algún poeta está a la altura de la fama de cualquiera, entonces, el último trabajo de cada individuo está a la altura de la fama de alguien.

Justificaciones:

| | | | |
|---|---|---|---|
| 1. *RR*. | 2, 3. *R*→ en 1. | 4. *IA*∃ y *A*∃ en 2. | 5, 6. *A*∧ en 4. |
| 7. *IR*∀ y *R*∀ en 3. | 8, 9. *IR*∃ y *R*∃ en 7. | 10, 11. *IA*∀ y *A*∀ en 6. | 12. *ABM*. |

Observar que las marcas de las fórmulas asociadas a las hojas determinan una interpretación $I = (D, v)$ que

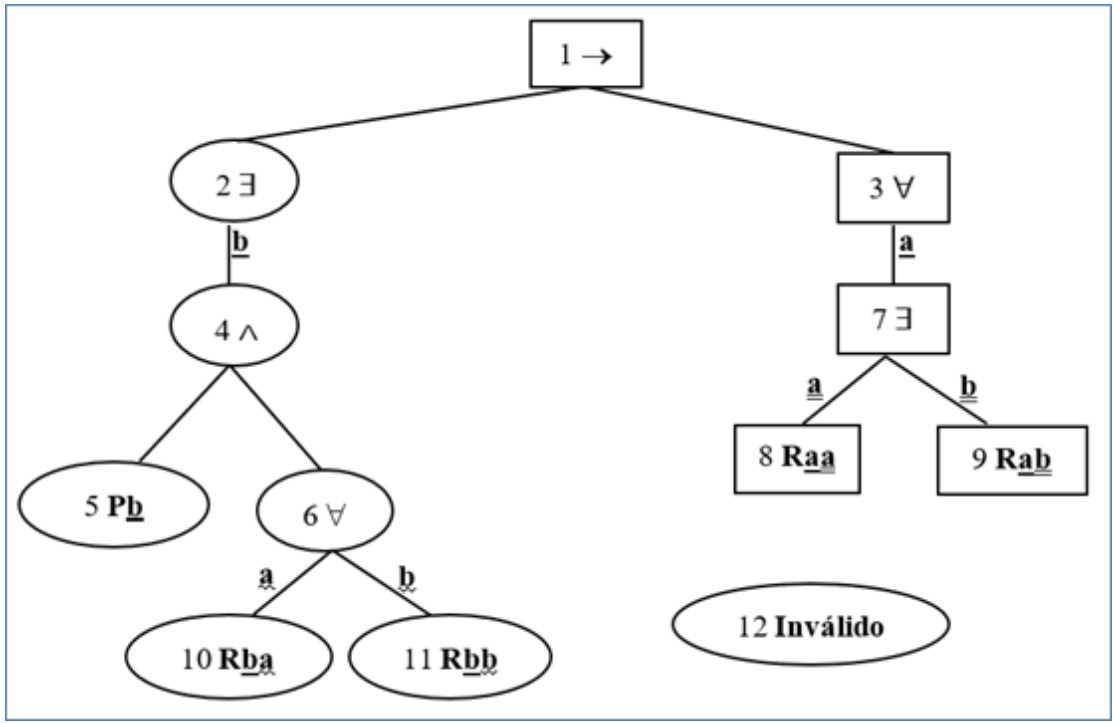

Figura 6: Árbol de forzamiento para la fórmula $\exists x(Px \wedge \forall yRxy) \rightarrow \forall x\exists yRxy$. Fuente: Elaboración Propia.

refuta a la fórmula analizada: $D = \{a, b\}, v(P) = \{b\}, v(R) = \{(b, a), (b, b)\}$.

La interpretación que refuta la validez de este argumento consta de dos individuos $a$ y $b$: Solo para fines ilustrativos, $a$ y $b$ se interpretan como Arturo y Bernardo, tales que Bernardo es un poeta, y su último trabajo está a la altura de la fama de Arturo y a la altura de su propia fama, mientras que Arturo no es un poeta, además, su último trabajo no está a la altura de la fama de Bernardo, ni a la altura de su propia fama.

Observar que no existe un modelo con un solo individuo que refute la fórmula, puesto que, necesariamente $a = b$, ya que del paso 8 se tiene $(a, a) \notin R$ pero del paso 11 se tiene $(b, b) \in R$.

Respecto al modelo que refuta, como el existencial del paso 7 es rechazado, el alcance de este existencial no puede ser satisfecho por ninguno de los individuos del modelo refutador, por esta razón deben salir dos ramas de este existencial, una por cada individuo del modelo; de igual manera, como el universal del paso 6 es aceptado, el alcance de este universal debe ser satisfecho por todos los individuos del modelo refutador, es decir, deben salir dos ramas de este universal.

## 5.2. Ilustración 2

En la figura 7 se muestra un árbol de forzamiento mal marcado para la fórmula $A$-válida $\forall x(\exists yPy \rightarrow \sim \exists y \sim Ryx) \rightarrow \sim \exists x(Px \wedge \sim \forall yRxy)$.

1. RR.
2, 3. R→ en 1.
4. R → en 3.
5. IA ∃ y A∧ en 4.
6,7 A∧ en 5.
8. A∼ en 7.
9. IR∀ y R∀ en 8.
10. IA∀ y A∀ en 2.
11. I∃ y IA en 6.
12. Aa∃ en 11.
13. AiA → en 12 y 10.
14. A∼ en 13.
15. IR∃ y R∃ en 14.
16. R∼ en 15.
17. DM en 9 y 16.

Partiendo de las justificaciones que garantizan la validez del argumento, y utilizando las técnicas de argumentación deductiva presentadas en Sierra (2010), se puede construir una prueba por reducción al absurdo en el lenguaje natural: Supóngase que [2] ∀$x$(∃$yPy$ →∼ ∃$y$ ∼ $Ryx$). Además, si se supone que [3, 4] ∃$x$($Px$∧ ∼ ∀$yRxy$), llamando, $a$, a tal individuo, resulta que [5, 6] $Pa$ y [7, 8] ∼ ∀$yRay$, llamando, $b$, al individuo que no cumple, se sigue que [9] ∼ $Rab$. Por otro lado, del supuesto inicial, en particular ocurre que [10] ∃$yPy$ ∼ ∃$y$ ∼ $Ryb$, y como se tiene [6, 11] $Pa$, es decir [12] ∃$yPy$, se infiere que [13, 14] ∼ ∃$y$ ∼ $Ryb$, y en particular, se puede afirmar que [15, 16] $Rab$, lo cual contradice el resultado previo [9] ∼ $Rab$, en consecuencia, se tiene lo contrario del segundo supuesto, es decir [17] ∼ ∃$x$($Px$ ∧ ∼ ∀$yRxy$). Se ha probado de esta manera que [18] ∀$x$(∃$yPy$ →∼ ∃$y$ ∼ $Ryx$) →∼ ∃$x$($Px$∧ ∼ ∀$yRxy$).

Observación. En el paso 10: $IA$∀ y $A$∀ en 2, también debe instanciarse para la constante a, generándose una nueva rama del árbol. Como el objetivo es encontrar una contradicción entre las marcas de los nodos, lo cual se logra entre los pasos 9 y 16 (lo cual significa que el argumento es válido, por lo tanto, no existe un modelo que refuta su validez), resultando que la instanciación para la constante $a$ es irrelevante, por esta razón tal instanciación puede ser omitida (lográndose de esta manera una apreciación visual más sencilla). Un razonamiento similar aplica en el paso 15: $IR$∃ y $R$∃ en 14. Cuando el modelo que refuta realmente existe, para su construcción si es necesaria la instanciación de todas las constantes.

### 5.3. Ilustración 3

En la figura 8 se muestra un árbol de forzamiento con la raíz marcada con 1 para la fórmula $A$-válida ∀$x$(($Px$ ∧ $Qb$) ∧ ∃$yRxy$) → ∀$x$ ∼ ($Px$ →∼ $Qb$).

Justificaciones:

1. $OA$ (opción de aceptación.
2. $IA$∀ y $A$∀ en 1.
3, 4. $A$∧ en 2.
5. IA ∃ y $A$∧ en 4.
6, 7 A∧ en 3.
8. $IA$ en 7.
9. $Aa$ ∼ en 8.
10. $I$∀ y $IA$ en 6.
11. $AiRd$ → en 10 y 9.
12. $Ra$ ∼ en 11.
13. $Aa$∀ en 12.
14. $OAi − Ad$ → en 1, 13.
15. Raíz marcada con 1.

Observar en el paso 13, que para aplicar la regla $Aa$∀, se requiere que la variable $x$ sea independiente en el paso 12, lo cual es cierto ya que, $x$ es independiente del supuesto del paso 1 ($x$ no ocurre libre en 1), y además, $x$ es independiente de la constante $b$ en 11 ($b$ no es un testigo introducido por las reglas $IRa$∀ y $R$∀ o $IA$∃ y $A$∃).

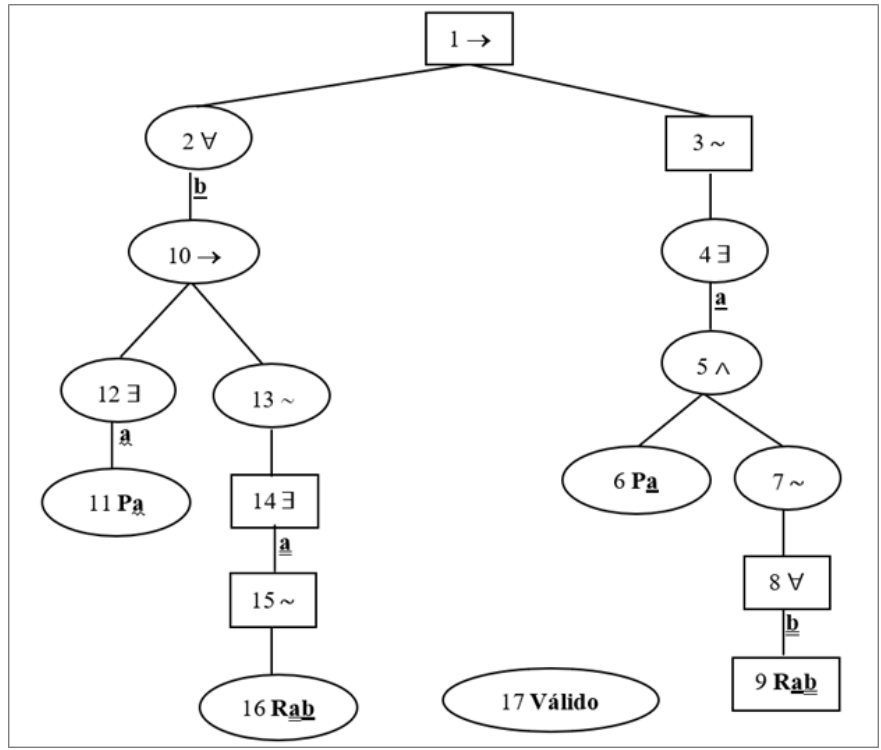

Figura 7: Árbol de forzamiento para la fórmula $\forall x(\exists y Py \to \sim \exists y \sim Ryx) \to \sim \exists x(Px \wedge \sim \forall y Rxy)$. Fuente: Elaboración Propia.

En el paso 2. *IA*$\forall$ y *A*$\forall$ en 1. Como no existen variables ni constantes, según la regla *IA*$\forall$, debe ser introducida alguna de ellas, en este caso (según la observación 3 de la sección 4.1) se debe introducir una variable para poder aplicar la regla *Aa*$\forall$ en el futuro paso 13.

## 5.4. Ilustración 4

En la figura 9 se muestra un árbol de forzamiento bien marcado para la fórmula *A*-inválida [$\forall x(x \in H \to x \in B) \wedge \exists x(\sim (x \in B) \wedge x \in A)$] $\to \forall x(x \in H \to \sim (x \in A))$, la cual corresponde a la siguiente fórmula de la teoría de conjuntos: $(H \subseteq B \wedge B^c \cap A \neq 0/) \to H \subseteq A^c$, o al siguiente argumento de la lógica tradicional: Todos los humanos son bípedos, algunos animales no son bípedos, por lo que ningún humano es un animal.

| | | | |
|---|---|---|---|
| 1. RR. | 2, 3. R$\to$ en 1. | 4. IR $\forall$ y R$\forall$ en 3. | 5, 6. R $\to$ en 4. |
| 7. R$\sim$ en 6. | 8, 9. A$\wedge$ en 2. | 10. IA$\exists$ y A$\exists$ en 9. | 11, 12. A$\wedge$ en 10. |
| 13. A $\sim$ en 11. | 14. IA $\forall$ y A $\forall$ en 8. | 15. IA en 5. | 16. AiA $\sim$ en 15 y 14. |
| 17. IA$\forall$ y A$\forall$ en 8. | 18. IR en 13. | 19. RdA $\to$ en 18 y 17. | 20. ABM . |

Observar que las marcas de las fórmulas asociadas a las hojas determinan una interpretación $I = (D, v)$, y en consecuencia los conjuntos que refutan a la fórmula analizada: $D = \{d, e\}, v(H) = Humanos = \{d\}$, $v(B) = $ Bípedos $= \{d\}$, $v(A) = $ Animales $= \{d, e\}$. De lo anterior se tiene que $B^c = \{e\}$, $A^C = \{\}$, de donde claramente $H \subseteq B$ y $B^c \cap A = \{e\} \neq \emptyset$, pero no se cumple que $\{d\} = H \cap A^c = \emptyset$.

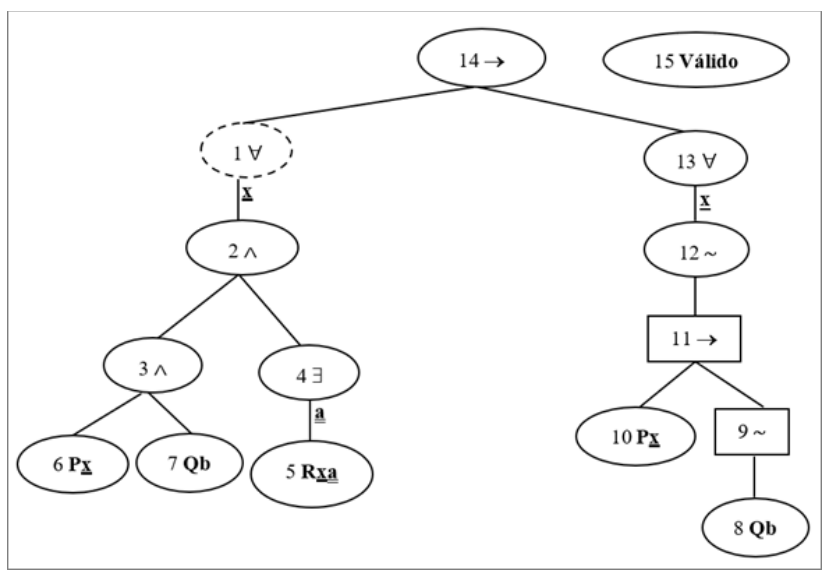

Figura 8: Árbol de forzamiento para la fórmula $\forall x((Px \land Qb) \land \exists y Rxy) \to \forall x \sim (Px \to \sim Qb)$. Fuente: Elaboración Propia.

Respecto al número de individuos, no existe un modelo con un solo individuo que refute la fórmula, puesto que necesariamente $d \neq e$, ya que del paso 5 se tiene $d \in H$ pero del paso 19 se tiene $e \notin H$.

Respecto a la construcción del modelo que refuta, como el universal del paso 8 es aceptado, el alcance de este universal debe ser satisfecho por todos los individuos del modelo refutador, por esta razón deben salir dos ramas de este universal, una por cada individuo del modelo.

### 5.5. Ilustración 5

En la figura 10 se muestra un árbol de forzamiento con la raíz marcada con 1 para la fórmula $A$-válida $\exists y \forall x Pyx \to \sim \exists x \forall y \sim Pyx$, la cual corresponde a la formalización del siguiente argumento: Si alguien acepta a todos los individuos, entonces, no es cierto que alguien sea rechazado por todos los individuos.

| | | | |
|---|---|---|---|
| 1. OA (Opción de aceptación) | 2. IA$\exists$ y A$\exists$ en 1. | 3. IA $\forall$ y A$\forall$ en 1. | 4. IA en 3. |
| 5. Aa $\sim$ en 4. | 6. RA$\forall$ en 5. | 7. Ra $\exists$ y en 6. | 8. Ra $\sim$ en 7. |
| 9. OAi-Ad $\to$ en 1, 8. | 10. Raiz marcada con 1. | | |

Observar en el paso 7, que para aplicar la regla $Ra\exists$, se requiere que la variable $x$ sea independiente en el paso 6, lo cual es cierto ya que, $x$ es independiente del supuesto del paso 1 ($x$ no ocurre libre en 1), y además, $x$ es independiente de la constante $a$ en 6 (razón 1: $a$ es un testigo introducido por las reglas $IA\exists$ y $A\exists$ en el paso 2, pero la variable $x$ no ocurre libre en el paso 2, razón 2: el testigo $a$ no figura en el paso 6).

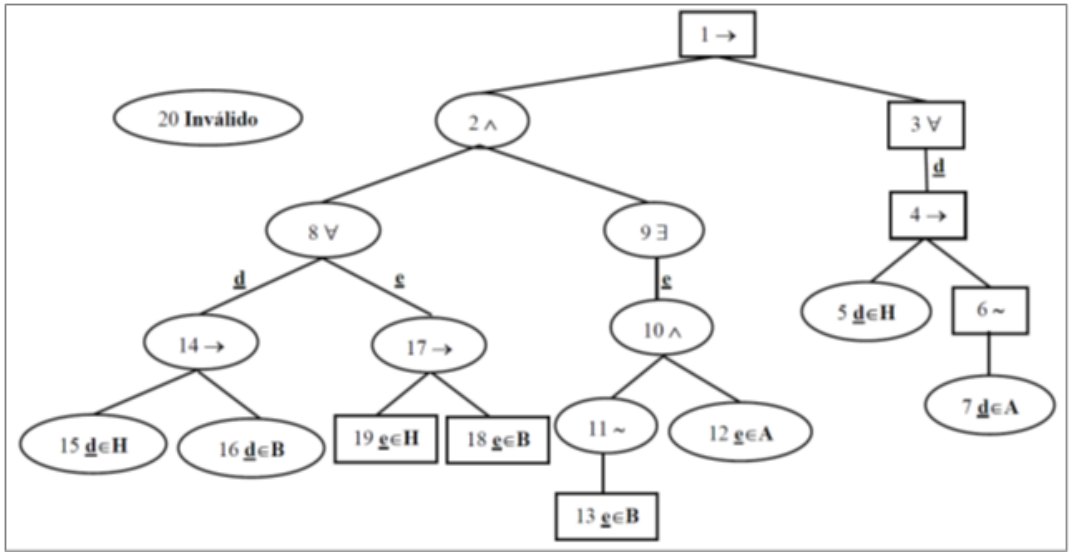

Figura 9: Árbol de forzamiento para $[\forall x(x \in H \to x \in B) \land \exists x(\sim x \in B \land x \in A)] \to \forall x(x \in H \to \sim x \in A)$. Fuente: Elaboración Propia.

Partiendo de las justificaciones que garantizan la validez del argumento, y utilizando las técnicas de argumentación deductiva presentadas en Sierra (2010), se puede construir una prueba directa en el lenguaje natural: Supóngase que, [1] alguien acepta a todos los individuos, llamando $a$ al tal alguien, resulta que [2] [$a$ acepta a todos los individuos], es decir, [3, 4] [$a$ acepta al individuo $x$], siendo $x$ un individuo cualquiera, por lo que, [5] [es falso que $a$ rechace al individuo $x$], en consecuencia, [6] [no es cierto que todos rechacen al individuo $x$], por lo tanto, [7, 8] [no es cierto que alguien sea rechazado por todos los individuos]. Se concluye finalmente que, [9] [si alguien acepta a todos los individuos, entonces, no es cierto que alguien sea rechazado por todos los individuos].

### 5.6. Ilustración 6

En la figura 11 se muestra un árbol de forzamiento bien marcado para la fórmula $A$-inválida $\forall x \exists y Pxy \to \exists y \forall x Pxy$, la cual corresponde a la formalización del siguiente argumento: Si todos donan algo a alguien, entonces, hay alguien a quien todos donan algo.

| | | | |
|---|---|---|---|
| 1. RR. | 2, 3. R $\to$ en 1. | 4. A $\forall$ en 3. | 5. A$\exists$. |
| 6. A $\forall$ en 2. | 7, 8. R$\exists$ en 3. | 9. OR en 8. | 10. OA en 6. |
| 11. OR-Ad en 7. | 12. ABM. | | |

Observar que en el paso 4, se tiene una rama para un individuo $a$ (el cual, según la regla $IA\forall$, existe porque los modelos no pueden ser vacíos). En el paso 6, se debe generar otra rama para el individuo nuevo *b* (el cual fue generado en el paso 5). Del paso 3 se deben generar dos ramas, una para cada individuo. En el paso 9, si se aplica la regla $R\forall$, se debe generar un individuo nuevo, lo cual implica que se deben generar nuevas ramas en los pasos 2 y 3, y el proceso se repite generando un número infinito de individuos; por esta razón, se intenta construir un modelo con dos individuos que refute la validez del argumento. En el paso 9 no se puede elegir el individuo *a*, ya que entraría en contradicción con el paso 5, por lo que se elige *b* como una opción de rechazo, en el paso 10 no se puede elegir el individuo *b*, ya que entraría en contradicción con el paso 9, por lo que se elige *a* como una opción de aceptación, en el paso 11 no se puede elegir el individuo *b*, ya que entraría en contradicción con el paso 10, por lo que se elige *a* como una opción de rechazo. A partir de las marcas de las hojas se tiene el modelo refutador: *a* dona algo a *b*, *b* dona algo a *a*, pero *a* no dona algo a *a*, ni *b* dona algo a *b*, como se muestra en la tabla de la figura 11.

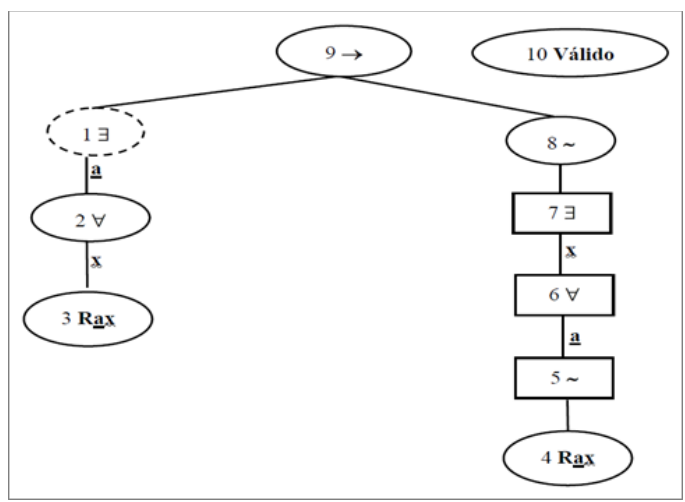

Figura 10: Árbol de forzamiento para $\exists y \forall x Pyx \rightarrow \sim \exists x \forall y \sim Pyx$. Fuente: Elaboración Propia.

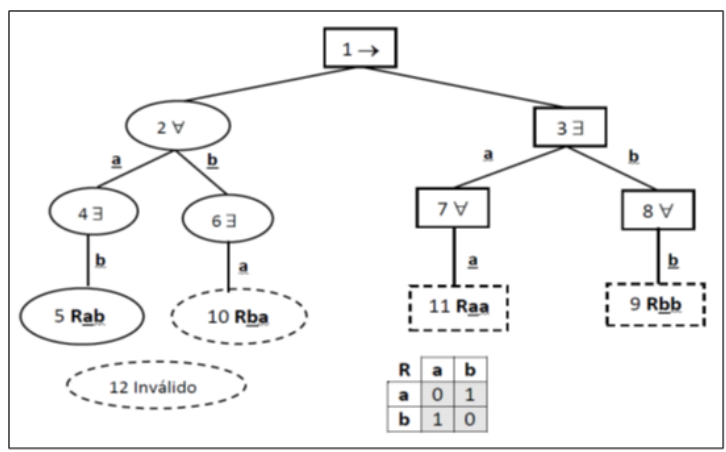

Figura 11: Árbol de forzamiento para $\forall x \exists y Pxy \rightarrow \exists y \forall x Pxy$. Fuente: Elaboración Propia.

- Observación 1. Este tipo de razonamiento puede ser recurrente cuando se tienen cuantificadores anidados, y para el caso de predicados no monádicos, no hay garantía de que el razonamiento funcione ya que en general la lógica de relaciones de primer orden no es decidible, Church (1936).

- Observación 2: Sin embargo, como la lógica de predicados monádicos si es decidible, Lowenheim (1915), entonces el razonamiento funciona, ya que se tiene el siguiente teorema: cuando en una fórmula $X$ se tienen $n$ predicados monádicos diferentes, si $X$ no puede ser refutada por un modelo con $2^n$ individuos, entonces $X$ no puede ser refutada por ningún modelo.

- Observación 3: Además, como la lógica de predicados diádicos con dos variables también es decidible, Börger *et al.* (1997), entonces el razonamiento funciona, ya que se tiene el siguiente teorema: cuando en una fórmula $X$ se tienen $n$ predicados monádicos y/o diádicos diferentes, si $X$ no puede ser refutada por un modelo con $2^{O(n)}$ individuos (lo cual puede lograrse en tiempo exponencial no determinista), entonces $X$ no puede ser refutada por ningún modelo. Para detalles consultar la proposición 8.1.4 de la referencia dada.

## 5.7. Ilustración 7

En la figura 12 se muestra un árbol de forzamiento (con un error en el paso 5) con el cual se prueba que la fórmula $\forall x \exists y Pxy \rightarrow \exists y \forall x Pxy$, es $A$-válida, la cual realmente es inválida (ver ilustración 6).

1. OA.    2. A $\forall$ en 1.    3. A $\exists$ en 2.    4. IA en 3.
5. Aa $\forall$ en 4.    6. Aa$\exists$ en 5.    7. OAi-Ad $\rightarrow$ en 1, 6.    8. RMI en 7.

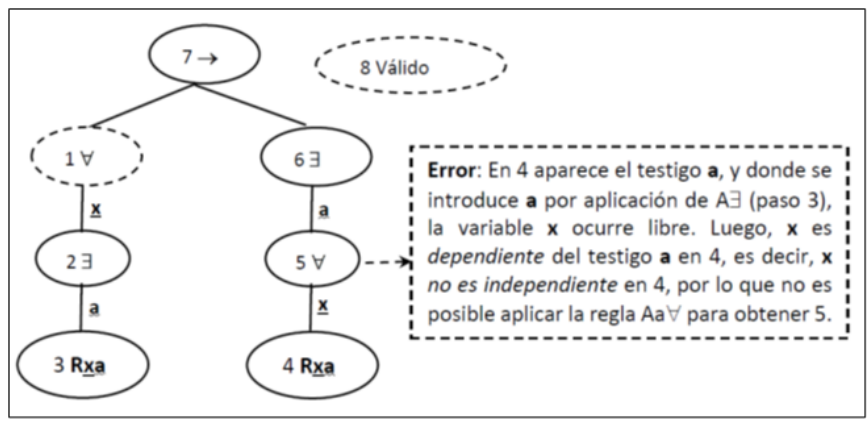

Figura 12: Árbol de forzamiento errado para $x \exists y Pxy \rightarrow y \forall x Pxy$. Fuente: Elaboración Propia.

# 6. MODELOS PARA LA LÓGICA DE PREDICADOS MONÁDICOS Y/O DIÁDICOS

Siguiendo las ideas de Henkin (1949) y adaptando la presentación hecha en Caicedo (1990), se tiene que una interpretación o modelo I de *LP*, es una pareja ordenada $I = (D, v)$, donde $D$ es un conjunto no vacío llamado el dominio del modelo, y $v$ es una valoración.

El lenguaje de *LP* incluye predicados monádicos y diádicos, variables y constantes.

- Si $P$ es un predicado monádico de *LP* entonces $v(P)$ es un subconjunto de $D$.
- Si $Q$ es un predicado diádico de *LP* entonces $v(P)$ es un subconjunto de $D^2$.
- Si $c$ es una constante de *LP* entonces $v(c)$ es un elemento fijo del dominio del modelo.
- Si $x$ es una variable de *LP* entonces $v(x)$ es un elemento arbitrario del dominio del modelo.

## 6.1. Reglas primitivas para la verdad de una fórmula en un modelo

Considerando el modelo $I = (D, v)$, sea $\ddot{a} = a_1, \ldots, a_n$ una secuencia de elementos del dominio $D$.

Si $P$ es un predicado diádico y $t_1, t_2$ son constantes o variables entonces se define: $I(P(t_1, t_2)[\ddot{a}]) = 1 \Leftrightarrow (v(t_1)[\ddot{a}], v(t_2)[]) \in v(P)$, donde, si $t_i$ es una variable entonces $v(t_i)[\ddot{a}] = a_i$, y si $t_i$ es una constante entonces $v(t_i)[\ddot{a}] = v(t_i)$.

Si $X$ y $Y$ son fórmulas cuyas variables libres se encuentran entre $x_1, \ldots, x_n$ entonces:

- *VI* $\sim$ . $I((\sim X)[\ddot{a}]) = 1 \leftrightarrow I(X[\ddot{a}]) = 0$.
- *VI* $\to$ . $I((X \to Y)[\ddot{a}]) = 0 \Leftrightarrow I(X[\ddot{a}]) = 1$ y $I(Y[\ddot{a}]) = 0$.
- *VI* $\wedge$ . $I((X \wedge Y)[\ddot{a}]) = 1 \Leftrightarrow I(X[\ddot{a}]) = I(Y[\ddot{a}]) = 1$.
- *VI* $\vee$ . $I((X \vee Y)[\ddot{a}]) = 0 \Leftrightarrow I(X[\ddot{a}]) = I(Y[\ddot{a}]) = 0$.
- *VI* $\leftrightarrow$ . $I((X \wedge Y)[\ddot{a}]) = 1 \Leftrightarrow I(X[\ddot{a}]) = I(Y[\ddot{a}])$.

Si las variables libres en la fórmula $F(x)$ son $x, x_1, \ldots, x_n$ entonces:

- *VI*$\exists$ 1. $I(\exists x F(x)[\ddot{a}]) = 1 \Rightarrow I(F(x)[b, \ddot{a}]) = 1$ para alguna $b$ nueva en $D$, $b$ llamada testigo.
- *VI*$\exists$ 2. $I(F(x)[b, \ddot{a}]) = 1$ para alguna $b$ en $D \Rightarrow I(\exists x F(x)[\ddot{a}]) = 1$.

- *VI∀* 1. $I(\forall x F(x)[\ddot{a}]) = 1 \Rightarrow I(F(x)[b, \ddot{a}]) = 1$ para toda $b$ en $D$.

- *VI∀* 2. $I(F(x)[b, \ddot{a}]) = 1$ para toda $b$ en $D$, con $b$ independiente en $F(x)[b, \ddot{a}] \Rightarrow I(\forall x F(x)[\ddot{a}]) = 1$.

En $F(x)[b, \ddot{a}]$, se dice que $b$ es independiente del supuesto $S(x)[b, \ddot{a}]$, si $F(x)[b, \ddot{a}]$ no se encuentra en el alcance del supuesto. Se dice que $b$ es independiente del testigo $e$ en $F(x, y)[b, e, \ddot{a}]$, si $b$ no ocurre donde el testigo $e$ fue introducido (donde aparece por primera vez por aplicación de *VI∃*1). Finalmente, se dice que en $F(x)[b, \ddot{a}]$, $b$ es independiente, si $b$ es independiente de todo testigo $e$ en $F(x, y)[b, e, \ddot{a}]$ y $b$ es independiente de todo supuesto $S(x)[b, \ddot{a}]$.

## 6.2. Validez de una fórmula en un modelo

- Si $X$ es una fórmula cuyas variables libres son $x_1, \ldots, x_n$ entonces:
    - $I(X) = 1 \Rightarrow I(X[\ddot{a}]) = 1$, para toda secuencia $\ddot{a}$ de elementos del dominio $D$.
    - $I(X) = 1$ se lee: la fórmula $X$ es verdadera en el modelo $I$.

- Si $X$ es una fórmula sin variables libres entonces:
    - $X$ es válida $\Rightarrow I(X) = 1$, para todo modelo o interpretación $I$.

- Si $X$ es una fórmula cuyas variables libres son $x_1, \ldots, x_n$ entonces:
    - $X$ es válida $\Rightarrow \forall x_1 \ldots \forall x_1 X$ es válida.

# 7. EQUIVALENCIA DE LAS PRESENTACIONES

Se define la Complejidad $C$, como una función la cual asigna a cada fórmula de *LP* un entero no negativo de la siguiente forma (donde $t_1, \ldots, t_n$ son variables o constantes):

- $C(P(t_1, \ldots, t_n)) = 0$, donde $P$ es un predicado $n$-ádico.

- $C(XkY) = 1 + $ máximo de $\{C(X), C(Y)\}$, donde $k \in \{\wedge, \vee, \rightarrow, \leftrightarrow\}$.

- $C(\sim X) = C(\forall_x X) = C(\exists_x X) = 1 + C(X)$.

Se define la Profundidad $P$, como una función la cual asigna a cada nodo de un árbol un entero no negativo de la siguiente forma (las señales, $s_1, \ldots, s_n$, pueden ser una variable, una constante o el espacio vacío, '_'):

- $P(P(s_1, \ldots, s_n)) = 0$, donde $P$ es un predicado $n$-ádico.

- $P(Ar[XkY]) = 1 + $ máximo de $\{P(Ar[X]), P(Ar[Y])\}$, donde $k \in \{\wedge, \vee, \rightarrow, \Leftrightarrow\}$.

- $P(Ar[\sim X]) = P(Ar[\forall_x X]) = P(Ar[\exists_x X]) = 1 + P(Ar[X])$.

### 7.1. Proposición 2. Unicidad de la extensión de una función de marca

Cada función de marca de hojas $m$, puede ser extendida de manera única, a una función de marca de nodos, $M$, de $N(X)$ en $\{0, 1\}$, haciendo $M(h) = m(h)$ si $h$ es una hoja, y aplicando las reglas de instanciación de espacios vacíos, junto con las reglas primitivas y derivadas para el forzamiento de marca, las cuales son presentadas en las secciones 4.1 a 4.5.

Prueba: Se debe probar que:

- La extensión $M$ de $m$ se aplique a todas las fórmulas.

- La asignación de $M$ a cada fórmula sea única.

- No existe otra extensión $M'$ de $m$, la cual se aplica a todas las fórmulas.

- Parte a.
  Por la definición, la extensión $M$ de $m$ se aplica a todas las fórmulas.

- Parte b.
  Se prueba por inducción sobre la profundidad $P$ del árbol de las fórmulas.
  Paso base. Profundidad 0, significa que el nodo es una hoja, en este caso $M$ coincide con la función $m$, por lo que se satisface la unicidad de asignación.

- Paso inductivo. Sea $P(Ar[X]) = L$, con $L > 0$. Como hipótesis inductiva se tienen: $P(Ar[Y]) < L$ implica la asignación de $M$ a $Y$ es única, $P(Ar[Z]) < L$ implica la asignación de $M$ a $Z$ es única, $P(Ar[Wa]) < L$ implica la asignación de $M$ a $Wa$ es única.

  - Caso 1.
    $X = \sim Z$. Supóngase que $M$ asigna 1 a $Ar[\sim Z]$, y que $M$ asigna 0 a $Ar[\sim Z]$. Como $M$ asigna 1 a $Ar[\sim Z]$, por la regla $A \sim$ se deriva que $M$ asigna 0 a $Ar[Z]$, además, como $M$ asigna 0 a $Ar[Z]$, por la regla $R \sim$ se deriva que $M$ asigna 1 a $Ar[Z]$, resultando que la asignación de $M$ a $Ar[Z]$ no es única, lo cual contradice la hipótesis inductiva. Por lo tanto, la asignación de $M$ a $\sim Z$ es única.

  - Caso 2.
    $X = Y \wedge Z$. Supóngase que $M$ asigna 1 a $Ar[Y \wedge Z]$, y que $M$ asigna 0 a $Ar[Y \wedge Z]$. Como $M$ asigna 1 a $Ar[Y \vee Z]$, por la regla $A\wedge$ se obtiene que $M$ asigna 1 a $Ar[Y]$ y $M$ asigna 1 a $Ar[Z]$, se tiene entonces que $M$ asigna 0 a $Ar[Y \wedge Z]$ y $M$ asigna 1 a $Ar[Y]$, aplicando la regla $RiA\wedge$ se infiere que $M$ asigna 0 a $Ar[Z]$, pero $M$ asigna 1 a $Ar[Z]$, lo cual contradice la hipótesis inductiva. Por lo tanto, la asignación de $M$ a $Y \wedge Z$ es única.

  - Caso 3.
    $X = Y \vee Z$. Supóngase que $M$ asigna 1 a $Ar[Y \vee Z]$, y que $M$ asigna 0 a $Ar[Y \vee Z]$. Como $M$ asigna 0 a $Ar[Y \vee Z]$, por la regla $R\vee$ se deriva que $M$ asigna 0 a $Ar[Y]$ y $M$ asigna 0 a $Ar[Z]$, se tiene

entonces que *M* asigna 1 a *Ar*[*Y* ∨ *Z*] y *M* asigna 0 a Ar[Y], aplicando la regla *RiA*∨ se infiere que *M* asigna 1 a *Ar*[*Z*], pero *M* asigna 0 a Ar[Z], lo cual contradice la hipótesis inductiva. Por lo tanto, la asignación de *M* a *Y* ∨ *Z* es única.

- Caso 4.

  *X* = *Y* → *Z*. Supóngase que *M* asigna 1 a *Ar*[*Y* → *Z*], y que *M* asigna 0 a *Ar*[*Y* → *Z*]. Como *M* asigna 0 a *Ar*[*Y* → *Z*], por la regla *R* → se deriva que *M* asigna 1 a *Ar*[*Y*] y *M* asigna 0 a *Ar*[*Z*], se tiene entonces que *M* asigna 1 a *Ar*[*Y* → *Z*] y *M* asigna 1 a *Ar*[*Y*], aplicando la regla *AiA* → se infiere que *M* asigna 1 a *Ar*[*Z*], pero *M* asigna 0 a *Ar*[*Z*], lo cual contradice la hipótesis inductiva. Por lo tanto, la asignación de *M* a *Y* → *Z* es única.

- Caso 5.

  *X* = *Y* ↔ *Z*. Supóngase que *M* asigna 1 a *Ar*[*Y* ↔ *Z*], y que *M* asigna 0 a *Ar*[*Y* ↔ *Z*]. Como *M* asigna 1 a *Ar*[*Y* ↔ *Z*], se consideran dos subcasos, *M* asigna 1 a *Ar*[*Y*] o *M* asigna 0 a *Ar*[*Y*].

  - Subcaso 1.

    *M* asigna 1 a *Ar*[*Y*]. Como *M* asigna 1 a *Ar*[*Y* ↔ *Z*], por la regla *AiA* ↔ se deduce que *M* asigna 1 a *Ar*[*Z*], como, además, *M* asigna 0 a *Ar*[*Y* ↔ *Z*], utilizando la regla *AdR* ↔ se concluye que *M* asigna 0 a *Ar*[*Y*], lo cual contradice la hipótesis inductiva.

  - Subcaso 2.

    *M* asigna 0 a *Ar*[*Y*]. Como *M* asigna 0 a *Ar*[*Y* ↔ *Z*], por la regla *RiR* ↔ se deduce que *M* asigna 1 a *Ar*[*Z*], como además, *M* asigna 1 a *Ar*[*Y* ↔ *Z*], utilizando la regla *AdA* ↔ se concluye que *M* asigna 1 a *Ar*[*Y*], lo cual contradice la hipótesis inductiva. Por lo tanto, la asignación de *M* a *Y* ↔ *Z* es única.

- Caso 6.

  *X* = ∀*xWx*. Supóngase que *M* asigna 1 a *Ar*[∀*xWx*], y que *M* asigna 0 a *Ar*[∀*xWx*]. Como *M* asigna 0 a *Ar*[∀*xWx*], utilizando las reglas *IR*∀ y *R*∀ se sigue la existencia de un testigo *a*, tal que *M* asigna 0 a *Ar*[*Wa*], pero como *M* asigna 1 a *Ar*[∀*xWx*], aplicando las reglas *IA*∀ y *A*∀ se deriva que *M* asigna 1 a *Ar*[*Wa*], lo cual contradice la hipótesis inductiva. Por lo tanto, la asignación de *M* a ∀*xWx* es única.

- Caso 7.

  *X* = ∃*xWx*. Supóngase que *M* asigna 1 a *Ar*[∃*xWx*], y que *M* asigna 0 a *Ar*[∃*xWx*]. Como *M* asigna 1 a *Ar*[∃*xWx*], utilizando las reglas *IA*∃ y *A*∃ se sigue la existencia de un testigo *a*, tal que *M* asigna 1 a *Ar*[*Wa*], pero como *M* asigna 0 a *Ar*[∃*xWx*], aplicando las reglas *IR*∃ y *R*∃ se deriva que *M* asigna 0 a *Ar*[*Wa*], lo cual contradice la hipótesis inductiva. Por lo tanto, la asignación de *M* a ∃*xWx* es única.

Por el principio de inducción matemática, se ha probado que la asignación de *M* a cada fórmula es única.

- Parte c.

Supóngase que existe otra extensión $M'$ de $m$, la cual es una función que se aplica a todas las fórmulas. Si $M \neq M'$, entonces existe al menos una fórmula $F$, tal que $M(Ar[F]) \neq M'(Ar[F])$. Entre estas fórmulas, el principio del buen orden garantiza la existencia de al menos una fórmula de profundidad mínima, sea $X$ una de estas fórmulas, por lo que $P(Ar[X]) = L$ y es mínima, es decir, para cada fórmula $T$, si $P(Ar[T]) < L$ entonces $M(Ar[T]) = M'(Ar[T])$. Además, como $M$ y $M'$ son ambas extensiones de $m$, entonces $Ar[X]$ no puede ser una hoja, por lo que $X$ debe ser una fórmula compuesta.

- Caso 1.

  $X = \sim Z$. Supóngase que $M(Ar[\sim Z]) = 1$ y que $M'(Ar[\sim Z]) = 0$. Como $M(Ar[\sim Z]) = 1$, aplicando la regla $A \sim$ se infiere que $M(Ar[Z]) = 0$, además, como $M'(Ar[\sim Z]) = 0$, utilizando la regla $R \sim$ se infiere que $M'(Ar[Z]) = 1$, por lo que $M(Ar[Z]) \neq M'(Ar[Z])$, lo cual no es posible puesto que $P(Ar[Z]) < L$.

- Caso 2.

  $X = Y \wedge Z$. Supóngase que $M(Ar[Y \wedge Z]) = 1$ y que $M'(Ar[Y \wedge Z]) = 0$. Como $M(Ar[Y \wedge Z]) = 1$, aplicando la regla $A\wedge$ se infiere que $M(Ar[Y]) = 1$ y $M(Ar[Z]) = 1$, además $P(Ar[Y]) < L$, de donde $M'(Ar[Y]) = 1$, se tiene entonces que $M'(Ar[Y \wedge Z]) = 0$ y $M'(Ar[Y]) = 1$ utilizando la regla $AiR\wedge$ se infiere que $M'(Ar[Z]) = 0$, por lo que $M(Ar[Z]) \neq M'(Ar[Z])$, lo cual no es posible puesto que $P(Ar[Z]) < L$.

- Caso 3.

  $X = Y \vee Z$. Supóngase que $M(Ar[Y \vee Z]) = 1$ y que $M'(Ar[Y \vee Z]) = 0$. Como $M(Ar[Y \vee Z]) = 0$, aplicando la regla $R\vee$ se infiere que $M(Ar[Y]) = 0$ y $M(Ar[Z]) = 0$, además $P(Ar[Y]) < L$, de donde $M'(Ar[Y]) = 0$, se tiene entonces que $M'(Ar[Y \vee Z]) = 1$ y $M'(Ar[Y]) = 0$ utilizando la regla $RiA\vee$ se infiere que $M'(Ar[Z]) = 1$, por lo que $M(Ar[Z]) \neq M'(Ar[Z])$, lo cual no es posible puesto que $P(Ar[Z]) < L$.

- Caso 4.

  $X = Y \rightarrow Z$. Supóngase que $M(Ar[Y \rightarrow Z]) = 1$ y que $M'(Ar[Y \rightarrow Z]) = 0$. Como $M(Ar[Y\ rightarrow Z]) = 0$, aplicando la regla $R \rightarrow$ se infiere que $M(Ar[Y]) = 1$ y $M(Ar[Z]) = 0$, además $P(Ar[Y]) < L$, de donde $M'(Ar[Y]) = 1$, se tiene entonces que $M'(Ar[Y \rightarrow Z]) = 1$ y $M'(Ar[Y]) = 1$ utilizando la regla $AiA \rightarrow$ se infiere que $M'(Ar[Z]) = 1$, por lo que $M(Ar[Z]) \neq M'(Ar[Z])$, lo cual no es posible puesto que $P(Ar[Z]) < L$.

- Caso 5.

  $X = Y \leftrightarrow Z$. Supóngase que $M(Ar[Y \leftrightarrow Z]) = 1$ y que $M'(Ar[Y \leftrightarrow Z]) = 0$. Supóngase que $M(Ar[Y]) = 1$, y como $M(Ar[Y \leftrightarrow Z]) = 0$, aplicando la regla $AiR \leftrightarrow$ se infiere que $M(Ar[Z]) = 0$, además $P(Ar[Z]) < L$, de donde $M'(Ar[Z]) = 0$, se tiene entonces que $M'(Ar[Y \leftrightarrow Z]) = 1$ y $M'(Ar[Z]) = 0$ utilizando la regla $RdA \leftrightarrow$ se infiere que $M'(Ar[Y]) = 0$, por lo que $M(Ar[Y]) \neq M'(Ar[Y])$, lo cual no es posible puesto que $P(Ar[Y]) < L$, en consecuencia, $M(Ar[Y]) = 0$, y como $M(Ar[Y \leftrightarrow Z]) = 0$, aplicando la regla $RiR \leftrightarrow$ se infiere que $M(Ar[Z]) = 1$, además $P(Ar[Z]) < L$, de donde $M'(Ar[Z]) = 1$, se tiene entonces que $M'(Ar[Y \leftrightarrow Z]) = 1$ y $M'(Ar[Z]) = 1$ utilizando la regla $AdA \leftrightarrow$ se infiere que $M'(Ar[Y]) = 1$, por lo que $M(Ar[Y])\neq M'(Ar[Y])$, lo cual no es posible puesto que $P(Ar[Y]) < L$.

- Caso 6.

    $X = \forall xZx$. Supóngase que $M(Ar[\forall xZx]) = 1$ y que $M'(Ar[\forall xZx]) = 0$. Como $M'(Ar[\forall xZx]) = 0$, utilizando las reglas $IR\forall$ y $R\forall$ se deriva la existencia de un testigo $a$, tal que $M'(Ar[Za]) = 0$, y además $P(Ar[Za]) < L$, por lo que $M(Ar[Za]) = 0$, como $M(Ar[\forall xZx]) = 1$, por las reglas $IA\forall$ y $A\forall$ se concluye que $M(Ar[Za]) = 1$, por lo que $M(Ar[Za]) \neq M'(Ar[Za])$, lo cual no es posible puesto que $P(Ar[Za]) < L$.

- Caso 7. $X = \exists xZx$. Supóngase que $M(Ar[\exists xZx]) = 1$ y que $M'(Ar[\exists xZx]) = 0$. Como $M(Ar[\exists xZx]) = 1$, utilizando las reglas $IA\exists$ y $A\exists$ se deriva la existencia de un testigo $a$, tal que $M(Ar[Za]) = 1$, y además $P(Ar[Za]) < L$, por lo que $M'(Ar[Za]) = 1$, como $M'(Ar[\exists xZx]) = 0$, por las reglas $IR\exists$ y $R\exists$ se concluye que $M'(Ar[Za]) = 0$, por lo que $M(Ar[Za])\neq M'(Ar[Za])$, lo cual no es posible puesto que $P(Ar[Za]) < L$. Por lo tanto, no existe otra extensión $M'$ de $m$, la cual se aplica a todas las fórmulas.

## 7.2. Proposición 3. Modelo asociado a una función de marca

Para cada fórmula $X$ de LP y para cada función de marca $m$, existe un modelo $I_m$, tal que, $M(R[X]) = 1 \Leftrightarrow I_m(X) = 1$.

Prueba: Sea $X$ una fórmula de LP y sea $m$ una función de marcas para $X$. Se define el modelo $I_m = (D_m, v_m)$ de la siguiente forma:

- Si la constante $c$ figura en el árbol de la fórmula $X$ entonces $c \in D_m$ ($c$ puede ser una constante original de $X$, o un testigo generado por las reglas $IA\exists$ o $IR\forall$).

- Si $c$ es una constante entonces $v_m(c) = c$.

- Si $x$ es una variable entonces $v_m(x) = x$ un elemento arbitrario de $D_m$.

- Si $P$ es un predicado monádico y $z_1$ es una variable o constante entonces $I_m(Pz_1) = 1 \Leftrightarrow z_1 \in v_m(P) \Leftrightarrow m(Pz_1) = 1$.

- Si $R$ es un predicado diádico y $z_1, z_2$ son variables o constantes entonces $I_m(Rz_1z_2) = 1 \Leftrightarrow (z_1, z_2) \in v_m(R) \Leftrightarrow m(Rz_1z_2) = 1$.

La verdad de las fórmulas en el modelo $I_m$ se define mediante las reglas primitivas de la sección modelos para la lógica de predicados. Se tiene entonces que $I_m$ es un modelo de LP.

Para probar que $M(R[X]) = 1 \Leftrightarrow I_m(X) = 1$, se procede por inducción sobre la complejidad de la fórmula $X$. Paso base: Supóngase que la $C(X) = 0$, esto significa que $X = Pz_1$ o $X = Qz_1z_2$ donde $P$ es un predicado monádico, $Q$ es un predicado diádico y $z_1, z_2$ son variables o constantes.

- Caso 1.

  Al ser $X$ atómica, se tiene que $M(R[Pz_1]) = 1 \Leftrightarrow M(Pz_1) = 1$, además en el modelo $I_m$, $M(Pz_1) = 1 \Leftrightarrow z_1 \in v_m(P)$, por la definición de $vm$ resulta $z_1 \in v_m(P) \Leftrightarrow v_m(z_1) \in v_m(P)$, como para $c$ constante, $x$ variable y $a\ddot{}$ una secuencia arbitraria de elementos de $D_m$ se tiene que $v_m(c[a\ddot{}]) = v_m(c) = c$ y $v_m(x[a\ddot{}]) = v_m(x) = x$, entonces resulta que $v_m(z_1) \in v_m(P) \Leftrightarrow v_m(z_1[a\ddot{}]) \in v_m(P)$, por la definición de verdad se tiene $v_m(z_1[a\ddot{}]) \in v_m(P) \Leftrightarrow I_m(Pz_1[a\ddot{}]) = 1$, y al ser $a\ddot{}$ una secuencia arbitraria de elementos del dominio se sabe que $I_m(Pz_1[a\ddot{}]) = 1 \Leftrightarrow I_m(Pz_1) = 1$. Se ha probado de esta manera que $M(R[Pz_1]) = 1 \Leftrightarrow I_m(Pz_1) = 1$, es decir, $M(R[X]) = 1 \Leftrightarrow I_m(X) = 1$.

- Caso 2.

  Al ser $X$ atómica, se tiene que $M(R[Qz_1z_2]) = 1 \Leftrightarrow M(Qz_1z_2) = 1$, además en el modelo $I_m$, $M(Qz_1z_2) = 1 \Leftrightarrow (z_1, z_2) \in v_m(Q)$, por la definición de $v_m$ resulta $(z_1, z_2) \in v_m(Q) \Leftrightarrow (v_m(z_1), v_m(z_2)) \in v_m(Q)$, como para $c$ constante, $x$ variable y $a\ddot{}$ una secuencia arbitraria de elementos de $D_m$ se tiene que $v_m(c[a\ddot{}]) = v_m(c) = c$ y $v_m(x[a\ddot{}]) = v_m(x) = x$, entonces resulta que $(v_m(z_1), v_m(z_2)) \in v_m(Q) \Leftrightarrow (vm(z_1[a\ddot{}]), v_m(z_1[a\ddot{}])) \in v_m(Q)$, por la definición de verdad se tiene $(v_m(z_1[a\ddot{}]), v_m(z_1[a\ddot{}])) \in v_m(Q) \Leftrightarrow I_m(Pz_1z_2[a\ddot{}]) = 1$, y al ser $a\ddot{}$ una secuencia arbitraria de elementos del dominio se sabe que $I_m(Qz_1z_2[a\ddot{}]) = 1 \Leftrightarrow I_m(Qz_1z_2) = 1$. Se ha probado de esta manera que $M(R[Qz_1z_2]) = v1 \Leftrightarrow I_m(Qz_1z_2) = 1$, es decir, $M(R[X]) = 1 \Leftrightarrow I_m(X) = 1$.

- Paso de Inducción:

  Supóngase que $C(X) \geq 1$. Al ser $C(X) \geq 1$, $X$ debe ser una fórmula compuesta, es decir, $X$ tiene una de las siguientes formas: $\forall x X(x), \exists x X(x)$ (en los casos en los cuales $X$ tiene una de las siguientes formas: $B \wedge C, B \vee C, B \to C, B \leftarrow C, \sim B$, se procede como en Sierra (2006)). Se analiza cada caso por separado.

  - Caso 1:

    Sea $X(x)$ una fórmula con variables libres $x, x_1, \ldots, x_n$, sea $Z = \forall x X(x)$. Se tiene que $R[Z] = \forall$, por lo que $M(R[Z]) = 1 \Leftrightarrow M(\forall) = 1$, pero por las reglas $IA\forall$ y $A\forall$ se tiene $M(\forall) = 1 \Rightarrow M(a\forall) = 1$, y como además $a\forall = R[X(x)[b]]$ para toda constante o variable $b$, resulta que $M(R[Z]) = 1 \Rightarrow M(R[X(x)[b]]) = 1$. Utilizando la hipótesis inductiva y el hecho que las constantes forman parte del dominio de $I_m$, se tiene que $M(R[Z]) = 1 \Rightarrow I_m(X(x)[b]) = 1$, para toda $b$ en el dominio del modelo $I_m$. Por la definición de verdad resulta $M(R[Z]) = 1 \Rightarrow I_m(X(x)[b, a\ddot{}]) = 1$, para toda $b$ en $D_m$ y para la secuencia arbitraria $a\ddot{}$ en $D_m$, y en consecuencia, al no utilizarse las reglas $R\forall$ ni $A\exists$ y no tener supuestos, para las $b$ en $D_m$, con $b$ independiente en $X(x)[b, a\ddot{}]$. Utilizando la regla $VI\forall 2$ se obtiene $M(R[X]) = 1 \Rightarrow I_m(\forall x X(x)[a\ddot{}]) = 1$, para la secuencia arbitraria $a\ddot{}$ en

$D_m$. Aplicando la definición de verdad se concluye $M(R[Z]) = 1 \Rightarrow I_m(\forall xX(x)) = 1$, es decir, $M(R[Z]) = 1 \Rightarrow I_m(Z) = 1$.

Para probar la recíproca, se sabe que $I_m(Z) = 1 \Leftrightarrow I_m(\forall xX(x)) = 1$, por la regla $VI\forall 1$ se obtiene $I_m(X) = 1 \Rightarrow I_m(X(x)[b]) = 1$ para toda $b$ en $D_m$, en consecuencia, al no utilizarse la regla $VI\exists$ y no tener supuestos, para una variable $b$ independiente. Utilizando la hipótesis inductiva, se tiene que $I_m(X) = 1 \Rightarrow M(R[X(x)[b]]) = 1$. Como $a\forall = R[X(x)[b]]$ resulta que $I_m(Z) = 1 \Rightarrow M(a\forall) = 1$ para una variable independiente en el espacio vacío. Finalmente, utilizando la regla $Aa\forall$, se obtiene $I_m(Z) = 1 \Rightarrow M(\forall) = 1$, es decir, $I_m(Z) = 1 \Rightarrow M(R[Z]) = 1$.

- Caso 2:

  Sea $X(x)$ una fórmula con variables libres $x, x_1, \ldots, x_n$, sea $Z = \exists xX(x)$. Se tiene que $R[Z] = \exists$, por lo que $M(R[Z]) = 1 \Leftarrow M(\exists) = 1$, pero por las reglas $IA\exists$ y $A\exists$ se tiene $M(\exists) = 1 \Rightarrow M(a\exists) = 1$, y como además $a\exists = R[X(x)[b]]$ para un nuevo testigo $b$, resulta que $M(R[Z]) = 1 \Rightarrow M(R[X(x)[b]]) = 1$. Utilizando la hipótesis inductiva y el hecho que los testigos forman parte del dominio de $Im$, se tiene que $M(R[Z]) = 1 \Rightarrow I_m(X(x)[b]) = 1$, para alguna $b$ en $D_m$. Por la definición de verdad resulta $M(R[Z]) = 1 \Rightarrow I_m(X(x)[b, a\ddot{}]) = 1$, para algún $b$ en $D_m$ y para la secuencia arbitraria $a\ddot{}$ en $D_m$. Utilizando la regla $VI\exists 2$ se obtiene $M(R[Z]) = 1 \Rightarrow I_m(\exists xX(x)[a\ddot{}]) = 1$, para toda secuencia $a\ddot{}$ en $D_m$. Aplicando la definición de verdad se concluye $M(R[Z]) = 1 \Rightarrow I_m(\exists xX(x)) = 1$, es decir, $M(R[Z]) = 1 \Rightarrow I_m(Z) = 1$.

  Para probar la recíproca, se sabe que $I_m(Z) = 1 \Leftrightarrow I_m(\exists xX(x)) = 1$, por la regla $VI\exists 1$ se obtiene $I_m(Z) = 1 \Rightarrow I_m(X(x)[b]) = 1$, para alguna $b$ nueva en $D_m$. Utilizando la hipótesis inductiva resulta $I_m(X) = 1 \Rightarrow M(R[X(x)[b]]) = 1$. Como $a = R[X(x)[b]]$ entonces $I_m(X) = 1 \Rightarrow M(a\exists) = 1$ para algún $b$ en el espacio vacío, aplicando la regla $Aa\exists$ resulta que $I_m(X) = 1 \Rightarrow M(\exists) = 1$, es decir, $Im(Z) = 1 \Rightarrow M(R[Z]) = 1$.

  Se tiene entonces que para todos los casos $M(R[Z]) = 1 \Leftrightarrow I_m(Z) = 1$, quedando así probado el paso de inducción. Por el principio de Inducción se concluye que: para toda fórmula $X$, $M(R[X]) = 1 \Leftrightarrow I_m(X) = 1$. Se ha probado entonces que para cada fórmula $X$ de *LP* y para cada función de marca $m$, existe un modelo $I_m$, tal que, $M(R[X]) = 1 \Leftrightarrow I_m(X) = 1$.

## 7.3. Proposición 4. Función de marca asociada a un modelo

Para cada fórmula $X$ de *LP* y para cada modelo $I$, existe una función de marca $m_I$, tal que, $MI_(R[X]) = 1 \Leftrightarrow I(X) = 1$.

Prueba: Sea $X$ una fórmula de LP y sea $I = (D, v)$ un modelo. Se define la función de marca mI del conjunto de hojas del árbol de $X$ en el conjunto $\{0, 1\}$ de la siguiente forma:

- Si P es un predicado monádico y $z_1$ es una variable o constantes (original o testigo) entonces $m_I(Pz_1) = 1 \Leftrightarrow v(z_1) \in v(P)$

- Si $Q$ es un predicado diádico y $z_1$ es una variable o constantes (original o testigo) entonces $m_I(Qz_1z_2) = 1 \Leftrightarrow (v(z_1), v(z_2)) \in v(Q)$

La función $m_I$ se extiende a una función $M_I$ del conjunto de nodos del árbol de $X$ en el conjunto 0, 1, mediante las reglas primitivas para el forzamiento de marcas presentadas en la sección marcando los nodos de un árbol. Se tiene entonces que $M_I$ es una función de marca de nodos. Para probar que $M_I(R[X]) = 1 \Leftrightarrow I(X) = 1$, se procede por inducción sobre la profundidad de Ar[X].

- Paso base: Supóngase que $P(Ar[X]) = 0$, esto significa que $X = Pz_1$ o $X = Qz_1z_2$ donde $P$ es un predicado monádico, $Q$ es un predicado diádico y $z_1$, $z_2$ son variables o constantes, y por lo tanto se tiene:

    • Caso 1.
    $M_I(R[X]) = M_I(R[Pz_1]) = M_I(Pz_1)$, y como se sabe que $M_I(Pz_1) = 1 \Leftrightarrow m_I(Pz_1) = 1$ y que $m_I(Pz_1) = 1 \Leftrightarrow v(z_1) \in v(P)$, se obtiene $M_I(R[X]) = 1 \Leftrightarrow v(z_1) \in v(P)$. Como para $c$ constante, $x$ variable y $a^{..}$ una secuencia arbitraria de elementos de $D$, se tiene que $v(c[a^{..}]) = v(c)$ y $v(x[a^{..}]) = v(x) = x$, se infiere $M_I(R[X]) = 1 \Leftrightarrow v(z_1)[a^{..}] \in v(P)$, lo cual significa que $M_I(R[X]) = 1 \Leftrightarrow I(Pz_1[a^{..}]) = 1$, utilizando la definición de verdad se concluye $M_I(R[X]) = 1 \Leftrightarrow I(Pz_1) = 1$, es decir, $M_I(R[X]) = 1 \Leftrightarrow I(X) = 1$.

    • Caso 2.
    $M_I(R[X]) = M_I(R[Pz_1z_2]) = M_I(Qz_1)$, y como se sabe que $M_I(Qz_1z_2) = 1 \Leftrightarrow m_I(Qz_1z_2) = 1$ y que $m_I(Qz_1z_2) = 1 \Leftrightarrow (v(z_1), v(z_2)) \in v(Q)$, se obtiene $M_I(R[X]) = 1 \Leftrightarrow (v(z_1), v(z_2))$. Como para $c$ constante, $x$ variable y $a^{..}$ una secuencia arbitraria de elementos de $D$, se tiene que $v(c[a^{..}]) = v(c)$ y $v(x[a^{..}]) = v(x) = x$, se infiere $M_I(R[X]) = 1 \Leftrightarrow (v(z_1)[a^{..}], v(z_2)[a^{..}]) \in v(Q)$, lo cual significa que $M_I(R[X]) = 1 \Leftrightarrow I(Qz_1z_2[a^{..}]) = 1$, utilizando la definición de verdad se concluye $M_I(R[X]) = 1 \Leftrightarrow I(Qz_1) = 1$, es decir, $M_I(R[X]) = 1 \Leftrightarrow I(X) = 1$.

- Paso de inducción: Supóngase que $P(Ar[X]) \geq 1$.
Al ser $P(Ar[X]) \geq 1$, $X$ debe ser una fórmula compuesta, es decir, $X$ tiene una de las siguientes formas: $\forall xX(x), \exists xX(x)$ (en los casos en los cuales $X$ tiene una de las siguientes formas: $B \wedge C, B \vee C, B \rightarrow C, B \leftarrow C, \sim B$, se procede como en Sierra (2006). Se analiza cada caso por separado.

    • Caso 1: Sea $X(x)$ una fórmula con variables libres $x, x_1, \cdots, x_n$, y sea $Z = \forall xX(x)$. Se tiene que $R[Z] = \forall$, por lo que $MI(R[Z]) = 1 \Leftrightarrow MI(\forall) = 1$, pero por las reglas $IA\forall$ y $A\forall \Rightarrow MI(a\forall) = 1$, y como además $a\forall = R[X(x)[b]]$ para toda constante o variable $b$, resulta que $M_I(R[Z]) = 1 \Rightarrow M_I(R[X(x)[b]]) = 1$. Utilizando la hipótesis inductiva y el hecho que las constantes forman parte del dominio del modelo $I$, se tiene que $M_I(R[Z]) = 1 \rightarrow I(X(x)[b]) = 1$, para toda b

en el dominio de $I$. Por la definición de verdad resulta $M_I(R[Z]) = 1 \Rightarrow I(X(x)[b, a\ddot{\,}]) = 1$, para toda $b$ en $D$ y la secuencia arbitraria $a\ddot{\,}$ en $D$, y al no utilizarse la reglas $R\forall$ ni $A\exists$ y no tener supuestos, para toda $b$ independiente en $X(x)[b, a\ddot{\,}]$. Utilizando la regla $VI\forall 2$ se obtiene $MI(R[Z]) = 1 \Rightarrow I(\forall x X(x)[a\ddot{\,}]) = 1$, para toda secuencia $a\ddot{\,}$ en $D$. Aplicando la definición de verdad se concluye $M_I(R[Z]) = 1 \Rightarrow I(\forall x X(x)) = 1$, es decir, $M_I(R[Z]) = 1 \Rightarrow I(Z) = 1$. Para probar la recíproca, se tiene que $I(Z) = 1 \Leftrightarrow, mI(\forall x X(x)) = 1$, por la definición de verdad resulta $I(Z) = 1 \Leftrightarrow I(\forall x X(x)[a\ddot{\,}]) = 1$, para cada secuencia $a\ddot{\,}$ de elementos del domino $D$. Aplicando la regla $VI\forall 1$ resulta $I(Z) = 1 \Rightarrow I(X(x)[b, a\ddot{\,}]) = 1$, para cada $b$ en $D$ y la secuencia arbitraria $a\ddot{\,}$ de elementos de $D$, y por la definición de verdad se obtiene $I(Z) = 1 \Rightarrow I(X(x)[b]) = 1$ para cada $b$ en $D$. Utilizando la hipótesis inductiva se infiere $I(Z) = 1 \Rightarrow M_I(R[X(x)[b]]) = 1$ para cada constante o variable $b$ en el espacio vacío, y como además $a\forall = R[X(x)[b]]$, entonces $I(X) = 1 \Rightarrow M_I(a\forall) = 1$ para toda constante o variable $b$ en el espacio vacío, y por lo tanto, al no aplicarse la regla $VI\exists$ y no tener supuestos, para una variable independiente en el espacio vacío. Aplicando la regla $Aa\forall$ se deduce que $I(X) = 1 \Rightarrow M_I(\forall) = 1$, es decir, que $I(Z) = 1 \Rightarrow M_I(Z) = 1$.

- Caso 2:

  Sea $X(x)$ una fórmula con variables libres $x, x1, \ldots, x_n$, sea $Z = \exists x X(x)$. Se tiene que $R[Z] = \exists$, por lo que $M_I(R[Z]) = 1 \Leftrightarrow M_I(\exists) = 1$, pero por las reglas $IA\exists$ y $A\exists$ se tiene $M_I(\exists) = 1 \Rightarrow M_I(a\exists) = 1$, y como además $a\exists = R[X(x)[b]]$ para un nuevo testigo $b$, y entonces resulta que $M_I(R[Z]) = 1 \Rightarrow M_I(R[X(x)[b]]) = 1$. Utilizando la hipótesis inductiva y el hecho que los testigos forman parte del dominio de $I$, se tiene que $M_I(R[Z]) = 1 \Rightarrow I(X(x)[b]) = 1$, para alguna $b$ en $D$. Por la definición de verdad resulta $M_I(R[Z]) = 1 \Rightarrow I(X(x)[b, a\ddot{\,}]) = 1$, para algún $b$ en $D$ y para la secuencia arbitraria $a\ddot{\,}$ en $D$. Utilizando la regla $VI\exists 2$ se obtiene $M_I(R[Z]) = 1 \Rightarrow I(\exists x X(x)[a\ddot{\,}]) = 1$, para toda secuencia $a\ddot{\,}$ en $D$. Aplicando la definición de verdad se concluye $M_I(R[Z]) = 1 \Rightarrow I(\exists x X(x)) = 1$, es decir, $M_I(R[Z]) = 1 \Rightarrow I(Z) = 1$.

  Para probar la recíproca, se sabe que $I(Z) = 1 \Leftrightarrow I(\exists x X(x)) = 1$. Por la definición de verdad resulta $I(Z) = 1 \Leftrightarrow I(\exists x X(x)[a\ddot{\,}]) = 1$ para toda secuencia $a\ddot{\,}$ en $D$. Aplicando la regla $VI\exists$ se obtiene $I(Z) = 1 \Rightarrow I(X(x)[b, a\ddot{\,}]) = 1$ para alguna $b$ nueva en $D$. Por la definición de verdad resulta $I(Z) = 1 \Rightarrow I(X(x)[b]) = 1$, y utilizando la hipótesis inductiva se infiere $I(Z) = 1 \Leftarrow M_I(R[X(x)[b]]) = 1$ para algún $b$. Como $a\exists = R[X(x)[b]]$ entonces $I(Z) = 1 \Rightarrow M_I(a\exists) = 1$ para algún $b$, y aplicando la regla $Aa\exists$ se deduce $I(Z) = 1 \Rightarrow M_I(\exists) = 1$, es decir, $I(Z) = 1 \Rightarrow M_I(Z) = 1$.

Se tiene entonces que para todos los casos $M_I(R[Z]) = 1 \Leftrightarrow I(Z) = 1$, quedando así probado el paso de inducción.

Por el principio de Inducción se concluye que: Para toda fórmula $X$, $M_I(R[X]) = 1 \Leftrightarrow I(X) = 1$. Se ha probado entonces que para cada fórmula $X$ de $LP$ y para cada modelo $I$, existe una función de marca $m_I$, tal que, $M_I(R[X]) = 1 \Leftrightarrow I(X) = 1.$

### 7.4. Proposición 5. Caracterización semántica de los árboles de forzamiento

- La fórmula $X$ es válida desde el punto de vista de los árboles si y solamente sí $X$ es válida desde el punto de vista de los modelos.

- La fórmula $X$ es válida desde el punto de vista de los árboles si y solamente sí $X$ es un teorema de la lógica de predicados monádicos y/o diádicos de primer orden (sin identidad ni funciones).

- Prueba parte a:
  Supóngase que $X$ no es válida desde el punto de vista de los árboles, entonces existe una función de marca $m$, tal que $M(R[X]) = 0.$ Se tiene entonces, por la proposición 3, que existe un modelo $I_m$, tal que $I_m(X) = 0$, y por lo tanto, $X$ no puede ser válida desde el punto de vista de los modelos.

  Supóngase ahora que $X$ no es válida desde el punto de vista de los modelos, entonces existe un modelo $I$, tal que $I(X) = 0$. Se tiene entonces, por la proposición 4, que existe una función de marca $m_I$, tal que $M_I(R[X]) = 0$, y por lo tanto, $X$ no puede ser válida desde el punto de vista de los árboles.

  Se concluye así que una $X$ es válida desde el punto de vista de los árboles si y solamente sí $X$ es válida desde el punto de vista de los modelos.

- Prueba parte b:
  De Caicedo (1990), se sabe que los modelos presentados caracterizan la lógica de predicados monádicos y/o diádicos, en consecuencia, una fórmula $X$ es válida desde el punto de vista de los árboles si y solamente sí $X$ es un teorema de la lógica de predicados monádicos y/o diádicos de primer orden.

## 8. CONCLUSIONES

Con los árboles de forzamiento semántico (forzamiento indirecto) para la lógica de predicados monádicos, y/o diádicos es posible determinar la validez de una fórmula de manera visual y relativamente mecánica, por ejemplo, utilizando un algoritmo para recorrer el árbol de la fórmula y en cada nodo buscando la aplicación de una regla para marcarlos (tal como se sugiere en el párrafo siguiente). Cuando una fórmula es inválida, es decir, cuando el árbol de la fórmula está bien marcado, entonces la lectura de las marcas de las hojas, tal como se hace en las ilustraciones 1, 4 y 6 (figuras 6, 9, 11), proporciona un modelo que refuta la validez de la fórmula.

Cuando la fórmula es válida, el árbol de forzamiento está mal marcado. Se tiene entonces que los árboles de forzamiento semántico (con el forzamiento indirecto) es un método refutacional, tal como lo son, los tableaux semánticos sistematizados por Smullyan. Las reglas de opciones para el condicional y la disyunción no son necesarias, pero facilitan la obtención del resultado final. Respecto a los árboles de forzamiento directo, estos no corresponden a un método refutacional.

Las marcas de los nodos de los árboles de forzamiento presentadas en las ilustraciones de la sección 5, se han realizado mediante la búsqueda intuitiva de nodos que puedan ser marcados. Pero este proceso puede ser sistematizado de manera más rigurosa, tal como se indica en el siguiente párrafo.

Definición. Un nodo está resuelto si se han aplicado todas las reglas que lo involucran. Sugerencia. Algoritmo para marcar los nodos de un árbol en el caso del forzamiento indirecto. Se hace el recorrido del árbol en preorden iniciando en la raíz, y se aplican los siguientes pasos hasta que termine el recorrido:

- Paso 1
  Marcar la raíz con 0.

- Paso 2
  Si el nodo no está resuelto, de ser posible, aplicar las reglas y marcar el nodo como resuelto.

- Paso 3
  Pasar al nodo siguiente.

- Paso 4
  Si el nodo no está resuelto, regresar al paso 2. Si el nodo está resuelto, regresar al paso 3.

  Una vez terminado el recorrido:

- Paso 4
  Si es posible determinar la validez o invalidez entonces el algoritmo termina.

- Paso 5
  Revisar los nodos marcados con 1 y asociados a un universal, deben salir tantas ramas como constantes, si hacen falta ramas éstas deben agregarse y marcarse con 1. Revisar los nodos marcados con 0 y asociados a un existencial, deben salir tantas ramas como constantes, si hacen falta ramas éstas deben agregarse y marcarse con 0.

- Paso 6
  Si es posible determinar la validez o invalidez entonces el algoritmo termina.

- Paso 7
  Si no es posible determinar la validez o invalidez, se toma una opción (de aceptación o rechazo) en alguno de los nodos, y se recorre de nuevo el árbol a partir del paso 2.

- Paso 8

  Si se obtiene una contradicción, se cambia el valor de la opción por la certeza del valor contrario, y se recorre de nuevo el árbol a partir del paso 2.

Si en el árbol figuran *n* predicados monádicos y ningún otro símbolo de predicado, entonces:

- Paso 9 Si después de generar $2^n$ individuos diferentes el algoritmo no termina, entonces no generar más *n* individuos nuevos, solo tomar opciones con los individuos existentes en las reglas generadoras de individuos.

De esta manera se garantiza la finalización del algoritmo, ya que, si el argumento no puede ser refutado con $2^n$ individuos diferentes, entonces no puede ser refutado, (Lowenheim, 1915).

Si en el árbol figuran *n* predicados monádicos y/o diádicos diferentes, solo dos variables y ningún otro símbolo de predicado, entonces:

- Paso 10

  Si después de generar $2^{O(n)}$ individuos diferentes el algoritmo no termina, entonces no generar más individuos nuevos, solo tomar opciones con los individuos existentes en las reglas generadoras de individuos. De esta manera se garantiza la finalización del algoritmo, ya que, si el argumento no puede ser refutado con $2^{O(n)}$ individuos diferentes entonces no puede ser refutado, (Börger *et al.*, 1997).

Para analizar la validez de un argumento relacionado con las operaciones entre conjuntos, se recurre a los diagramas de Venn presentados en Venn (1880), los círculos de Euler detallados en Hammer (1996) y las tablas de pertenencia utilizadas en Grimaldi (1998). Del procedimiento mostrado en la ilustración 4 (figura 9), se puede afirmar que los árboles de forzamiento semántico para predicados monádicos, y también para predicados diádicos con dos variables, proporcionan un método alternativo para el análisis de validez de argumentos del álgebra de conjuntos, y cuando el argumento es inválido, los árboles de forzamiento generan los conjuntos que sirven de contraejemplo. Aunque el argumento de álgebra de conjuntos mostrado en la ilustración 4, puede ser refutado mediante los árboles de forzamiento para operaciones entre conjuntos presentados en Sierra (2017), los árboles de forzamiento semántico presentados en este trabajo tienen un mayor alcance, al no limitarse al álgebra de conjuntos.

En Van Dalen (2004) se presenta un sistema de deducción natural para la lógica clásica de predicados. Las demostraciones en los sistemas de deducción natural son presentadas en forma de árbol, donde las premisas se presentan en las hojas y la conclusión en la raíz, y con frecuencia estos árboles también son llamados árboles de forzamiento. Sin embargo, los sistemas de deducción natural, al contrario de los árboles de forzamiento semántico, no proporcionan un método para refutar los argumentos inválidos, ya que no constituyen una herramienta semántica.

En Sierra (2010) se hace explícita la conexión directa que existe entre las reglas para el forzamiento de marcas (conectivos proposicionales) y las reglas de deducción natural. Esta conexión se extiende fácilmente a las reglas Árboles de Forzamiento Semántico para la Lógica de Predicados para los cuantificadores, por lo que, para el caso de argumentos válidos, el árbol de forzamiento semántico puede traducirse en una prueba directa en el lenguaje natural (como en la ilustración 5), o puede traducirse en una prueba por reducción al absurdo en el lenguaje natural (como en la ilustración 2).

Desde el punto de vista didáctico, el carácter intuitivo de las reglas para el forzamiento de marcas, hacen de los árboles de forzamiento semántico para la lógica de predicados, una herramienta de trabajo muy práctica.

## Agradecimientos

## Referencias